\documentclass[onecolumn,amsmath,showkeys,amssymb]{revtex4}

\usepackage{graphicx,color}

\usepackage{amssymb}
\usepackage{amsbsy}
\usepackage{amsthm,mathrsfs,amsopn}
\usepackage{amsmath}
\usepackage{mathtools}

\theoremstyle{plain} %% This is the default

\newtheorem{theorem}{Theorem}

\newtheorem{proposition}[theorem]{Proposition}

\newtheorem{remark}[theorem]{Remark}

\begin{document}

\title{Quantifying the degree of average contraction of Collatz orbits.}

\author{Timoteo Carletti$^{1}$, Duccio Fanelli$^{2}$}
\affiliation{
1. Department of Mathematics and Namur Center for Complex Systems - naXys, University of Namur, rempart de la Vierge 8, B 5000 Namur, Belgium\\
2. Dipartimento di Fisica e Astronomia, University of Florence,INFN and CSDC, Via Sansone 1, 50019 Sesto Fiorentino, Florence, Italy}

\begin{abstract} 
We here elaborate on a quantitative argument to support the validity of the Collatz conjecture, also known as the $(3x+1)$ or Syracuse conjecture. The analysis is structured as follows. First, three distinct fixed points are found for the third iterate of the Collatz map, which hence organise in a  period $3$ orbit of the original map.  These are  $1$, $2$ and $4$,  the elements which define the unique attracting cycle, as hypothesised by Collatz. To carry out the calculation we write the positive integers in modulo $8$ (\textit{mod8}),  obtain a closed  analytical form for the associated map and determine the transitions that yield contracting or expanding iterates in the original, infinite-dimensional, space of positive integers.  Then, we consider a Markov chain which runs on the reduced space of {\textit{mod8}} congruence classes of integers.  The transition probabilities of the Markov chain are computed from the deterministic map, by employing a measure that is invariant for the map itself. Working in this setting, we demonstrate that the stationary distribution sampled by the stochastic system induces a contracting behaviour for the orbits of the deterministic map on the original space of the positive integers. Sampling the equilibrium distribution on the congruence classes {\textit{mod$8^m$}} for any $m$, which amounts to arbitrarily reducing the degree of imposed coarse graining, returns an identical conclusion. 
\end{abstract}

\maketitle

\section{Introduction}

The Collatz conjecture is named after Lothar Collatz, who first proposed it in 1937~\cite{Lagarias}.  The conjecture is also known as the $(3x+1)$ conjecture, the Ulam conjecture (after Stanislaw Ulam), Kakutani's problem (after Shizuo Kakutani), the Thwaites conjecture (after Sir Bryan Thwaites), Hasse's algorithm (after Helmut Hasse), or the Syracuse problem \cite{Lagarias}.  This can be formulated as an  innocent problem of arithmetics. The beauty of the conjecture indeed emanates  from its apparent, tantalising, simplicity, which however hides formidable challenges for those who attempt to grasp its deep-rooted essence.  To date the conjecture is still unsolved, despite mathematicians strived hard to confer it the rigour of an indisputable proof. The problem can be stated as follow. Take any positive integer $n$. If $n$ is even, divide it by $2$ to get $n / 2$. If n is odd, multiply it by $3$ and add $1$ to obtain $3n + 1$. In formulae:
\begin{equation}
 \label{eq:themap}
 \forall n\in\mathbb{N}\quad T(n)=
 \begin{cases}
 \frac{n}{2}& \text{if $n$ is even}\\
3n+1& \text{if $n$ is odd.}
 \end{cases}
\end{equation}

Repeating the process iteratively, the map is believed to converge to a period--$3$ orbit formed by the triad $\{1,2,4\}$.  Equivalently, the conjecture states that the Collatz map will always reach $1$, no matter what integer number one starts with. Numerical experiments have confirmed the validity of the conjecture for extraordinarily large values of the starting integer $n$ \cite{oliveira}. 

In this paper, we provide a novel argument to support the validity of the Collatz conjecture. To anticipate our findings, we shall demonstrate that the third iterate $T^{\circ 3}$ of the Collatz map admits three fixed points, $1$, $2$ and $4$. These latter  elements define the, supposedly unique, attracting cycle conjectured by Collatz. The third iterate map is naturally defined on the \textit{mod8} congruence classes of positive integers. We thus quantify the factor of relative compression or expansion, as follows the application of $T^{\circ 3}$, on each of the eight congruence classes obtained under the \textit{mod8} operation.  In the second part of the paper, we show that  orbits are on average bound to asymptotically shrink in size so heading towards the deputed equilibrium. We will further enhance the resolution of the measure by reducing to an arbitrary extent the degree of imposed coarse graining, i.e. working on the 
the congruence classes {\textit{mod$8^m$}} for any chosen $m$. Working in this generalised setting, we will prove that the average Collatz dynamics is contracting, for  $m$ large as sought, namely shrinking the congruence classes arbitrarily close to the singletons corresponding to each integer.  

Markov processes based on congruence classes invariant under application of $T$ have been previously considered in \cite{Matthews1,Matthews2} and revisited by Lagarias in his comprehensive survey on the $(3x+1)$ problem~\cite{Lagarias}. It is however the combined usage of (i) the third iterate map  $T^{\circ 3}$, (ii) the representation of numbers in {\textit{mod8}}, (iii) the idea of employing a Markov chain constructed from $T^{\circ 3}$ via a suitably defined measure, that allows us to draw a rigorous bound for the contraction factor, which is not just heuristically guessed. 

To guide the reader through the text we shall hereafter provide a schematic outline of the main steps involved in the analysis, by making explicit reference to specific key results.

\begin{itemize} 
\item We will begin by defining the third iterate of the Collatz map hereby named $S=T^{\circ 3}$. 
\item We will determine the action of $S$ on integers expressed in $mod8$. We will obtain class-dependent, $\mathcal{B}(i,8) \, (i=0,\dots,7)$, expansion/contraction factors that exemplify the action of $S$, see Eqs.~(\ref{eq:the3map}). Working in this setting we will also show that $1,2,4$ are the only fixed points of the deterministic map $S$. The subsequent analysis is targeted to showing that the trajectories of $S$ are bound to converge on average to one of the above fixed points.
\item To this end we first introduce a finite states Markov chain which runs on the eight congruence classes $\mathcal{B}(i,8) \, (i=0,\dots,7)$. The transition probabilities are given by equation (\ref{eq:Qij}) and have been obtained using the $S$-invariant measure $\mu_{inv}$ on the classes $\mathcal{B}(j,8^m)$, $m\geq 1$ and $j=0,\dots,8^m-1$. 
\item The measure $\mu_{inv}$ is defined by equation (\ref{eq:mumeasure2}) (or equivalently equation (\ref{eq:meas_classes})).
The invariance of the measure under $S$ is proved in Theorem \ref{thm:invmeas}.
\item Since the transition probabilities are computed from the $S$-invariant measure, it is possible to draw conclusion on the iterates of $S$ (namely its restriction in $mod8$) by iterating forward the Markov process. This observation follows from a straightforward application of the Chapman-Kolmogorov equation, as discussed in Proposition  \ref{prop:invariance}. The explicit form of the stochastic matrix $Q^*$ that characterises the introduced Markov chain is given in Proposition \ref{prop:markov}.
The stochastic chain does not account  for the specificity of $1,2,4$, the equilibria of $S$.  It will hence allow us to elaborate on the out-of-equilibrium dynamics of $S$, prior  (possible) convergence to the asymptotic Collatz equilibrium. 
\item The stationary distribution of the Markov chain is computed and given by formula (\ref{eq:P_QS}). Recall that by iterating forward the Markov chain one can inspect the equilibrium dynamics of $S$, in its $mod8$ representation, see Proposition  \ref{prop:invariance}. 
\item By using the expansion/contraction factors associated to each of the classes $\mathcal{B}(i,8) \, (i=0,\dots,7)$ one can show that the deterministic trajectories are on average contracting. This is substantiated by formula  (\ref{eq:Qstella}).
\item The analysis is generalised by working on the congruence classes {\textit{mod$8^m$}}, for any $m$. By operating in this setting, we will proof that the average Collatz dynamics is contracting, for  arbitrarily large $m$, i.e. shrinking the size of the congruence classes as sought. Remarkably, the estimated upper bound for the contraction factor is shown to be independent on $m$.
\end{itemize}

\section{From the deterministic map to a stochastic framework.} 
The aim of this section is to introduce the tools used to derive our conclusion, through three intermediate steps, that we here outline for convenience. In subsection~\ref{ssec:3maps} we will define and analyse the third iterate of the deterministic Collatz map. Then, in subsections~\ref{ssec:markov},~\ref{ssec:invmeasN} and~\ref{ssec:transprob}  we will introduce a finite states Markov process, whose transition probabilities come from the aforementioned third iterate map. To construct the Markov chain we employ an invariant measure of the deterministic map.  In section~\ref{sec:wCollatz} we will study the stationary distribution of the proposed Markov chain and infer robust constraints for the Collatz (deterministic) dynamics. In particular, Collatz orbits are on average contracting.  Finally, we will expand on previous results by reformulating the dynamics  on the {\textit{mod$8^m$}}  classes, $\forall m$, and prove the average Collatz dynamics to be contracting at any level of imposed coarse graining.

\subsection{The third iterate of the Collatz map}
\label{ssec:3maps}

We begin by remarking that  $\{1,2,4\}$ is a period--$3$ orbit of the Collatz map. It is hence quite natural to operate with the third iterate of the original map. Working in this context, the elements of the Collatz cycle should emerge as distinct fixed points of the third iterate map. Furthermore, we progress in the analysis with a \textit{mod8} representation of the natural numbers, a choice that makes it possible to cast the sought map in a rather compact form. In practical terms, this amounts to organise the visited numbers in eight different congruence classes, $\mathcal{B}(i,8)$ for $i\in\{0,\dots, 7\}$, each containing the positive integers that yield the same remainder, after performing an Euclidean division by $8$. Mathematically, the 
 $\mathcal{B}(i,8)$ class is defined as:
\begin{equation}
 \label{eq:Ck}
\mathcal{B}(i,8):=\{n\in\mathbb{N}:\exists m\in\mathbb{N}\cup \{0\},\; n=i+8m\}\, .
\end{equation}

For the sake of simplicity,  we term $S$ the third iterate of the Collatz map, namely  $S:=T^{\circ 3}$. The following Proposition makes explicit the action of $S$ on a generic positive integer $n$. The outcome depends on the specific \textit{mod8} class, $n$ belongs to.
\begin{proposition}
 \label{prop:T3}
 Let $n$ be any positive integer belonging to $\mathcal{B}(i,8)$ for some $i\in\{0,\dots,7\}$, then the third iterate map $S$ is explicitly given by
 \begin{equation}
 \label{eq:the3map}
 \forall n\in\mathbb{N}\quad S(n)=
 \begin{cases}
\frac{n}{8}& \text{if $n\in \mathcal{B}(0,8)$}\\
\frac{6n+2}{8}& \text{if $n\in \mathcal{B}(1,8)$}\\
\frac{6n+4}{8}& \text{if $n\in \mathcal{B}(2,8)$}\\
\frac{36n+20}{8}& \text{if $n\in \mathcal{B}(3,8)$}\\
\frac{6n+8}{8}& \text{if $n\in \mathcal{B}(4,8)$}\\
\frac{6n+2}{8}& \text{if $n\in \mathcal{B}(5,8)$}\\
\frac{6n+4}{8}& \text{if $n\in \mathcal{B}(6,8)$}\\
\frac{36n+20}{8}& \text{if $n\in \mathcal{B}(7,8)$\, .}
 \end{cases}
\end{equation}
\end{proposition}

One can easily prove that the fixed points of $S$ are $1$, $2$ and $4$ corresponding to the $3$-cycle $\{1,2,4\}$ of the original Collatz map. This is achieved by solving the fixed point 
equation $S(n)=n$ in the interval of interest $n \ge 1$. It can be straightforwardly proved that {\it no additional fixed points exist} for the map $S(n)$.

Notice that Eq.~\eqref{eq:the3map} can be cast in a more compact form as follows
\begin{equation}
 \label{eq:the3mapbis}
S(n)=\frac{m_i n+r_i}{8} \text{ if }n\in\mathcal{B}(i,8) \, (i=0,\dots,7)\, ,
\end{equation}
where the integers $(m_i)_{0\leq i\leq 7}$ and $(r_i)_{0\leq i\leq 7}$ are given by
\begin{eqnarray}
 \label{eq:miri}
m_0&=&1, m_1=m_2=m_4=m_5=m_6=6 \text{ and }m_3=m_7=36\\
r_0&=&0, r_1=r_5=2, r_2=r_6=4, r_4=8 \text{ and }r_3=r_7=20\, .
\end{eqnarray}

As a side remark, we observe that $S$ is thus in the form of a generalised Collatz map~\cite{Matthews1}. In the following we will need the explicit value of $x_i=S(i)={(m_i i+r_i)}/{8}$ for $i=0,\dots,7$, that is
\begin{equation}
 \label{eq:xi}
x_0=0\, , x_1=1\, ,x_2=2\, ,x_3=16\, ,x_4=4\, ,x_5=4\, ,x_6=5 \text{ and }x_7=34\, .
\end{equation}

\subsection{A finite states Markov process}
\label{ssec:markov}

The (finite states) Markov process that we are going to introduce considers the aforementioned congruence classes, $\mathcal{B}(i,8)$, as a finite alphabet. The transition probabilities among different states follow the deterministic map $S(n)$, provided one works with a suitable probability space $(\mathbb{N},\mu_{inv})$, for some {S-}{\em invariant measure} $\mu_{inv}$, that we will introduce hereafter. More precisely, for any given pair of classes $\mathcal{B}(i,8)$ $(i=0,\dots,7)$ and $\mathcal{B}(j,8)$ $(j=0,\dots,7)$ the probability~\footnote{We shall hereafter label the rows and columns of $Q^*$ with indexes running from $0$ up to $7$, rather than from $1$ up to $8$, as it is customarily done.} $q^*_{ij}$ of being initially in $\mathcal{B}(i,8)$ and then land in $\mathcal{B}(j,8)$, that is the conditional probability $P\left[S(x)\in \mathcal{B}(j,8)| x\in \mathcal{B}(i,8)\right]$, is given by:
\begin{equation}
\label{eq:Qij}
q^*_{ij}:=\frac{\mu_{inv}\left[\mathcal{B}(i,8)\cap S^{-1}\mathcal{B}(j,8)\right]}{\mu_{inv}\left[\mathcal{B}(i,8)\right]}\, (i,j=0,\dots,7)\, .
\end{equation}

To compute the above transition probabilities, one needs to explicitly determine $S^{-1}\mathcal{B}(j,8)$, $j=0,\dots, 7$.  To gather this information we start with a preliminary remark: 

\begin{remark}[On the solution of congruence linear equations]
\label{rem:congeq}
Let us recall a basic fact of congruence linear equations; given integers $a$, $b$ and $n$, the equation
\begin{equation*}
ax\equiv b \quad \textit{mod n}\, ,
\end{equation*}
can be solved if and only if $d=gcd(a,n)$ (gcd stand for the greatest common divisor) is a divisor of $b$ (that is $b$ can be divided by $d$), in this case the number of distinct solutions is given by $d$.
\end{remark}

We are now in a position to prove the following result:

\begin{proposition}
 \label{prop:smenouno}
Let $j=0,\dots, 7$, then $S^{-1}\mathcal{B}(j,8)$ is the union of disjoint congruence classes \textit{mod64}, $\mathcal{B}(l_j,64)$, where the indexes $l_j$ depend on the {\textit{mod8}} congruence class $j$.

In explicit form:
\begin{eqnarray}
\label{eq:congclass864}
S^{-1}\mathcal{B}(0,8)&=&\mathcal{B}(0,64)\cup \mathcal{B}(10,64)\cup \mathcal{B}(42,64)\cup \mathcal{B}(3,64)\cup \mathcal{B}(19,64)\cup \mathcal{B}(35,64)\cup \mathcal{B}(51,64)\cup \mathcal{B}(20,64)\cup \mathcal{B}(52,64)\cup \notag\\
& &\cup \mathcal{B}(21,64)\cup\mathcal{B}(53,64)\notag\\
S^{-1}\mathcal{B}(1,8)&=&\mathcal{B} (1,64)\cup \mathcal{B}(33,64)\cup \mathcal{B}(22,64)\cup \mathcal{B}(54,64)\cup \mathcal{B}(8,64)\notag\\
S^{-1}\mathcal{B}(2,8)&=&\mathcal{B} (2,64)\cup \mathcal{B}(34,64)\cup \mathcal{B}(12,64)\cup \mathcal{B}(44,64)\cup \mathcal{B}(13,64)\cup \mathcal{B}(45,64)\cup \mathcal{B}(7,64)\cup \mathcal{B}(23,64)\cup \mathcal{B}(39,64)\cup \notag\\
& &\cup \mathcal{B}(55,64)\cup\mathcal{B}(16,64)\notag\\
S^{-1}\mathcal{B}(3,8)&=&\mathcal{B}(25,64)\cup \mathcal{B}(57,64)\cup \mathcal{B}(14,64)\cup \mathcal{B}(46,64)\cup \mathcal{B}(24,64)\notag\\
S^{-1}\mathcal{B}(4,8)&=&\mathcal{B}(26,64)\cup \mathcal{B}(58,64)\cup \mathcal{B}(11,64)\cup \mathcal{B}(27,64)\cup \mathcal{B}(43,64)\cup \mathcal{B}(59,64)\cup \mathcal{B} (4,64)\cup \mathcal{B}(36,64)\cup \mathcal{B}(5,64)\cup \notag\\
& &\cup \mathcal{B}(37,64)\cup\mathcal{B}(32,64)\notag\\
S^{-1}\mathcal{B}(5,8)&=&\mathcal{B}(17,64)\cup \mathcal{B}(49,64)\cup \mathcal{B}(6,64)\cup \mathcal{B}(38,64)\cup \mathcal{B}(40,64)\notag\\
S^{-1}\mathcal{B}(6,8)&=&\mathcal{B}(18,64)\cup \mathcal{B}(50,64)\cup \mathcal{B}(28,64)\cup \mathcal{B}(60,64)\cup \mathcal{B}(29,64)\cup \mathcal{B}(61,64)\cup \mathcal{B}(15,64)\cup \mathcal{B}(31,64)\cup \mathcal{B}(47,64)\cup \notag\\
& &\cup \mathcal{B}(63,64)\cup\mathcal{B}(48,64)\notag\\
S^{-1}\mathcal{B}(7,8)&=&\mathcal{B}(9,64)\cup \mathcal{B}(41,64)\cup \mathcal{B}(30,64)\cup \mathcal{B}(62,64)\cup \mathcal{B}(56,64)\, .
\end{eqnarray}
\end{proposition}

\proof
Let $n\in \mathcal{B}(l,64)$ for some $l=0,\dots,63$, that is it exists $k\in\mathbb{N}\cup\{0\}$ such that $n=l+64k$. Let $l\equiv i$ \textit{mod8}, namely $l=i+8h$ for some $i=0,\dots,7$ and $h=0,\dots,7$.

We can then evaluate $S(n)$ using Eq.~\eqref{eq:the3mapbis}:
\begin{equation*}
S(n)=\frac{m_i n+r_i}{8} =\frac{m_i (l+64k)+r_i}{8}=\frac{m_i l+r_i}{8}+{8m_ik}=\frac{m_i (i+8h)+r_i}{8}+{8m_ik}=x_i+m_ih+{8m_ik}\, ,
\end{equation*}
where we used the definition of $x_i$ in the rightmost step. Finally $S(n)\in\mathcal{B}(j,8)$ for some $j=0,\dots,7$, if and only if $S(n)\equiv j$ \textit{mod8}, that is
\begin{equation*}
x_i+m_ih\equiv j \quad \textit{mod8}\, ,
\end{equation*}
or equivalently 
\begin{equation}
\label{eq:mainceq}
m_ih\equiv j - x_i\quad \textit{mod8}\, .
\end{equation}

As stated in Remark~\ref{rem:congeq}, Eq.~\eqref{eq:mainceq} can be solved if and only if $d_i=gcd(m_i,8)$ is a divisor of $j - x_i$. 

Let us turn to compute $d_i$. From the definitions Eq.~\eqref{eq:miri} we readily get
\begin{equation}
\label{eq:di}
d_0=gcd(1,8)=1\, , d_1=d_2=d_4=d_5=d_6=gcd(6,8)=2 \text{ and }d_3=d_7=gcd(36,8)=4\, .
\end{equation}

Consider first the case $j=0$. Because $d_0=1$, Eq.~\eqref{eq:mainceq} has always one solution for $i=0$, given by $h=0$. Thus $l=i+8h=0$. For $i=1$ ($d_1=2$ and $x_1=1$) Eq.~\eqref{eq:mainceq} can be solved if and only if $j-1$ can be divided by $2$. Hence, it has no solution when $j=0$. For $i=2$ ($d_2=2$ and $x_2=2$) Eq.~\eqref{eq:mainceq} can be solved because $j-x_2=0-2$ is divisible by $2$. 
There are in particular $2$ solutions, $h=1$ and $h=5$, which return $l=i+8h=10$ and $l=42$. For $i=3$ ($d_3=4$ and $x_3=16$), Eq.~\eqref{eq:mainceq} can be solved if and only if $j-16$ can be divided by $4$. This holds true for $j=0$.  Four solutions are found which correspond to $h=0,2,4,6$, yielding $l=3,19,35,51$. For $i=4$ ($d_4=2$ and $x_4=4$), Eq.~\eqref{eq:mainceq} can be solved if and only if $j-4$ can be divided by $2$, and this is true if $j=0$. The two obtained solutions are $h=2$ and $h=6$ giving in turn $l=20$ and $l=52$. The same is true for $i=5$ ($d_5=2$ and $x_5=4$). In this case one obtains $l=21$ and $l=53$.  For $i=6$ ($d_6=2$ and $x_6=5$), Eq.~\eqref{eq:mainceq} can be solved if and only if $j-5$ can be divided by $2$ and this is impossible for $j=0$. Finally for $i=7$ ($d_7=4$ and $x_7=34$), Eq.~\eqref{eq:mainceq} can be solved if and only if $j-34$ can be divided by $4$ and this condition is not met for $j=0$. 

In conclusion we showed that  $S^{-1}\mathcal{B}(0,8)$ is the union of $\mathcal{B}(0,64)\cup \mathcal{B}(10,64)\cup \mathcal{B}(42,64)\cup \mathcal{B}(20,64)\cup \mathcal{B}(52,64)$ and $\mathcal{B}(3,64)\cup \mathcal{B}(19,64)\cup \mathcal{B}(35,64)\cup \mathcal{B}(51,64)\cup\mathcal{B}(21,64)\cup\mathcal{B}(53,64)$. The first set is composed by $5$ classes corresponding to even $l$ while the latter group is made by $6$ classes with odd $l$.

Let us now consider $j=1$. For $i=0$ there is only one solution $h=1$ and thus $l=8$. For $i=1$ there are two solutions, since $j-1$ is even. These are $h=0,4$, hence $l=1,33$. For $i=6$  two solutions are found because $j-5$ can be divided by $2$. These are $h=2,6$, yielding $l=22,54$. For the remaining cases $i=2,3,4,5,7$ no solutions are possible. Summing up, $S^{-1}\mathcal{B}(1,8)$ is the union of three classes with even $l$, $\mathcal{B}(22,64)\cup \mathcal{B}(54,64)\cup \mathcal{B}(8,64)$, and two classes with odd $l$, $\mathcal{B} (1,64)\cup \mathcal{B}(33,64)$.

The remaining cases can be handled similarly and are not discussed here in details.
\endproof

To go one step further and explicitly compute the transition probabilities given by Proposition~\ref{prop:smenouno} we introduce a {S-}{\em invariant (probability) measure}. As we will clarify in the following the measure is invariant under $S$ on the family of generalised congruence classes $B(j,8^m)$.

\subsection{A {S-}invariant probability measure}
\label{ssec:invmeasN}

For any fixed integer $m\geq 1$ we introduce a partition of $\mathbb{N}$ into disjoint congruence classes, $\mathcal{B}(i,8^m)$, $i\in\{0,\dots, 8^m-1\}$ and we will use them to define a ($1$-parameter family of) invariant measures $\mu^{(m)}_{inv}$. As we will clarify hereafter the introduced measure depends on the index $m$, which identifies the class of pertinence.

{Let $n$ be a positive integer and assume it belongs to $\mathcal{B}(i,8^m)$ for some $i$. Hence, $n=i+k\,8^m$ for some integer $k\geq 0$ and $i\in\{0,\dots, 8^m-1\}$. Then there exists a unique string $s_0,\dots,s_{m-1}\in\{0,\dots,7\}$ such that $i=s_{m-1}8^{m-1}+\dots+s_18+s_0$. Based on the above we define the measure  $\mu^{(m)}_{inv}$ of the integer $n$ to be:}
\begin{equation}
\label{eq:mumeasure2}
{\mu^{(m)}_{inv}(n)=\frac{1}{2^{k+1}}\frac{1}{8^{m-1}}\nu(s_0)\, ,}
\end{equation}
{where the factor $1/2^{k+1}$ is related to the considered partition \footnote{Observe that formally $k=k(m)$.} while the second factor quantifies the probability that $s_1,\dots,s_{m-1}$ take 
any of the symbolic entries $\{0,\dots,7\}$ in the expression for $i$.  These latter probabilities are assumed identical (i.e. equal $1/8$) for all symbols.}

{The quantities $\nu(s_0)$  are given by}
\begin{equation*}
{\nu(s_0)=\frac{1}{6}\text{ if $s_0=0,2,4,6$ and }\nu(s_0)=\frac{1}{12}\text{ if $s_0=1,3,5,7$.}}
\end{equation*}

{The measure of any set made by integers is the sum of the measures of the integers forming the set. The introduced measure is therefore additive. Based on the above we can straightforwardly measure any class $\mathcal{B}(i,8^m)$. Following  Eq.~\eqref{eq:mumeasure2} the measure of $\mathcal{B}(i,8^m)$ reads:}
\begin{equation*}
{\mu^{(m)}_{inv}(\mathcal{B}(i,8^m))=\sum_{n\in \mathcal{B}(i,8^m)}\mu^{(m)}_{inv}(n)=\sum_{k\geq 0 }\frac{1}{2^{k+1}}\frac{\nu(s_0)}{8^{m-1}}=\frac{\nu(s_0)}{8^{m-1}}\, ,}
\end{equation*}
{where $i=s_0+s_1 8+\dots+s_{m-1}8^{m-1}$ for $s_0,\dots,s_{m-1}\in\{0,\dots,7\}$, hence}
\begin{equation*}
{\mu^{(m)}_{inv}(\mathcal{B}(i,8^m))=\frac{1}{6}\frac{1}{8^{m-1}}\text{ if $s_0=0,2,4,6$ and }\mu^{(m)}_{inv}(\mathcal{B}(i,8^m))=\frac{1}{12}\frac{1}{8^{m-1}}\text{ if $s_0=1,3,5,7$.}}
\end{equation*}
{Clearly $\mathbb{N}=\bigcup_{s_0,\dots,s_{m-1}}\mathcal{B}(s_0+s_1 8+\dots+s_{m-1}8^{m-1},8^m)$ and thus}
\begin{equation}
\label{eq:mumN}
{\mu^{(m)}_{inv}(\mathbb{N})=\sum_{s_0,\dots,s_{m-1}} \mu^{(m)}_{inv}(\mathcal{B}(s_0+s_1 8+\dots+s_{m-1}8^{m-1},8^m))=\sum_{s_0,\dots,s_{m-1}}\nu(s_0)\frac{1}{8^{m-1}}=\sum_{s_0}\nu(s_0)\frac{1}{8^{m-1}}8^{m-1}=1\, ,}
\end{equation}
{that is $\mu^{(m)}_{inv}$ is a probability measure for any $m\geq 1$.}

\begin{remark}

Notice that the invariant measure $\mu^{(m)}_{inv}$ can be equivalently defined working at the  level of classes. One can in fact define:

\begin{equation}
\label{eq:meas_classes}
\mu^{(m)}_{inv}(\mathcal{B}(i,8^m)) = \frac{1}{8^{m-1}} \nu(i)
\end{equation}

where 

\begin{equation}
\label{meas_def}
\nu(i)=\frac{1}{6}\text{ if $i (mod8) = 0,2,4,6$ and }\nu(i)=\frac{1}{12}\text{ if $i (mod 8)=1,3,5,7$.}
\end{equation}

Clearly:

\begin{equation*}
 {\mu^{(m)}_{inv}(\mathbb{N})=\sum_{i=0}^{8^m-1} \mu^{(m)}_{inv}(\mathcal{B}(i,8^m))=\sum_{\sigma=0}^7\nu(\sigma)=1\, ,}
\end{equation*}

\end{remark}

One can prove that the above introduced measure is invariant under $S$ for all sets made by any finite intersection and union of the congruence classes defined above. This is established in the following Theorem.

\begin{theorem}
\label{thm:invmeas}
For any $m\geq 1$ and for all $j=0,\dots,8^m-1$ we have
\begin{equation}
\label{eq:muinvB}
\mu_{inv}\left[S^{-1}\mathcal{B}(j,8^m)\right]=\mu_{inv}\left[\mathcal{B}(j,8^m)\right]\, .
\end{equation}
\end{theorem}

The proof of the above theorem relies on the following Proposition
\begin{proposition}
\label{prop:AeAo}
For any $m\geq 1$ and for all $j=0,\dots,8^m-1$ we have
\begin{equation}
\label{eq:Sm1Bj}
S^{-1}\mathcal{B}(j,8^m)=A^{(m)}_e(j)\cup A^{(m)}_o(j)\, ,
\end{equation}
where $A_e(j)$ is the union of disjoint classes $\mathcal{B}(l,8^{m+1})$ with $l$ even and $A^{(m)}_o(j)$ is the union of disjoint classes $\mathcal{B}(l,8^{m+1})$ with $l$ odd. Moreover if $j$ is even then $A^{(m)}_e(j)$ contains five elements and $A^{(m)}_o(j)$ six elements, while if $j$ is odd then $A^{(m)}_e(j)$ contains three elements and $A^{(m)}_o(j)$ two elements.
\end{proposition}

\proof
Observe that Eq.~\eqref{eq:Sm1Bj} holds true for $m=1$ by Proposition~\ref{prop:smenouno}. Let us assume it is true for all $k\leq m-1$ and prove it for $k=m$, namely we have to prove that $S^{-1}\mathcal{B}(j,8^m)$ is the disjoint union of classes $\mathcal{B}(l,8^{m+1})$ and more precisely $\#A_e^{(m)}(j)=5$ if $j$ is even and $3$ if $j$ is odd, while $\#A_o^{(m)}(j)=6$ if $j$ is even and $2$ if $j$ is odd.

Let thus $n\in\mathcal{B}(l,8^{m+1})$ for some $l=0,\dots,8^{m+1}-1$, that is $n=l+8^{m+1}k$ for some positive integer $k$. Let $l\equiv i$ \textit{mod8}, $l=i+8h$, $h=0,\dots,8^m-1$. From the definition~\eqref{eq:the3map} we can compute
\begin{equation*}
S(n)=\frac{m_in+r_i}{8}=\frac{m_i l+r_i}{8}+m_i8^{m}k=x_i+m_i h+m_i8^{m}k\, ,
\end{equation*}
hence $S(n)\in\mathcal{B}(j,8^{m})$ for some $j=0,\dots,8^{m}-1$ if and only if $S(n)\equiv j$ \textit{mod $8^m$}, namely
\begin{equation*}
m_ih\equiv j-x_i \quad \textit{mod $8^m$}\, .
\end{equation*}

This equation can be solved if and only if $d_i^{(m)}=gcd(m_i,8^m)$ is a divisor of $j-x_i$. Let us observe firstly that for all $m\geq 1$, $d_i^{(m)}=d_i$ where $d_i$ have been defined above Eq.~\eqref{eq:di}.

Let us now determine $j_1=0,\dots, 8^{m-1}$ such that $j\equiv j_1$ \text{mod $8^{m-1}$}. Observing that $j-x_i=j_1-x_i+q8^{m-1}$, for some integer $q$, we can conclude that $d_i$ divides $j-x_i$ if and only if it divides $j_1-x_i$, hence using the induction hypothesis we can conclude that $m_ih\equiv j-x_i \quad \textit{mod $8^m$}$ has solutions if and only if the same equation for $m-1$ has solutions, and moreover they have the same number of solutions.

This implies that $\#A_e^{(m+1)}(j)=\#A_e^{(m)}(j_1)$ and $\#A_e^{(m+1)}(j)=\#A_e^{(m)}(j_1)$, because the  \text{mod $8^{m-1}$} operation doesn't change the parity of $j$. The claim is then proved using the results of Proposition~\ref{prop:smenouno}.
\endproof

We can now prove Theorem~\ref{thm:invmeas}.
\proof
Let $m\geq 1$ and $j=0,\dots,8^{m}-1$, then thanks to Proposition~\ref{prop:AeAo} we have
\begin{equation*}
\mu_{inv}\left[S^{-1}\mathcal{B}(j,8^m)\right]=\mu_{inv}\left[A^{(m)}_e(j)\right]+\mu_{inv}\left[A^{(m)}_o(j)\right]\, ,
\end{equation*}
being the classes disjoint. Observe that all classes $\mathcal{B}(l,8^{m+1})\in A^{(m)}_e(j)$ have the same measure given by $1/6\times1/8^m$, and the same is true for classes in $A^{(m)}_o(j)$, with measure $1/12\times1/8^m$, thus
\begin{equation*}
\mu_{inv}\left[S^{-1}\mathcal{B}(j,8^m)\right]=\frac{1}{8^m}\frac{1}{6}\#A^{(m)}_e(j)+\frac{1}{8^m}\frac{1}{12}\#A^{(m)}_o(j)\, .
\end{equation*}

If $j$ is even we get:
\begin{equation*}
\mu_{inv}\left[S^{-1}\mathcal{B}(j,8^m)\right]=\frac{1}{8^m}\frac{1}{6}5+\frac{1}{8^m}\frac{1}{12}6=\frac{1}{8^{m-1}}\frac{1}{6}=\mu_{inv}\left[\mathcal{B}(j,8^m)\right]\, ,
\end{equation*}
If $j$ is odd we get:
\begin{equation*}
\mu_{inv}\left[S^{-1}\mathcal{B}(j,8^m)\right]=\frac{1}{8^m}\frac{1}{6}3+\frac{1}{8^m}\frac{1}{12}2=\frac{1}{8^{m-1}}\frac{1}{12}=\mu_{inv}\left[\mathcal{B}(j,8^m)\right]\, .
\end{equation*}

\endproof

\begin{remark}
Given a positive integer $m^*$ and the measure $\mu_{inv}^{(m^*)}$, all measures $\mu_{inv}^{(m)}$, $\forall m<m^*$ are push-forward measures of $\mu_{inv}^{(m^*)}$: $\mu_{inv}^{(m^*-1)}(\cdot):=\mu_{inv}^{(m^*)} (S^{-1}(\cdot))$.  In fact 
$\mu_{inv}^{(m^*-1)}(\mathcal{B}(i,8^{m^*-1})) = \mu_{inv}^{(m^*)} (S^{-1} \mathcal{B}(i,8^{m^*-1})) = \mu_{inv}^{(m^{*})} (\bigcup_j \mathcal{B}(j,8^{(m^{*})})) = \sum_j \mu_{inv}^{(m^*)} (\mathcal{B}(j,8^{m^*}))$
where the last sum extends over the classes $8^{m^*}$ that constitute the image of $\mathcal{B}(i,8^{m^*-1})$ via application of $S^{-1} (\cdot)$. Based on Proposition \ref{prop:AeAo}
it readily follows that $\mu_{inv}^{(m^*-1)}(\mathcal{B}(i,8^{m^*-1})) = \nu(i)/8^{m^*-1}$. The reasoning can be further iterated to all $m<m^*$, so recovering the measure of the classes as introduced above.  
\end{remark}

%\begin{remark}
%\ins{E' un riassunto di un mail di Vassilis, lo metto per non dimenticare, vediamo se vale la pena metterlo. To summarise for any given integer $m \geq 1$ the classes $\mathcal{B}(i, 8^m)$, $i = 0, \cdots, 8^m - 1$, form a partition of $\mathbb{N}$. Let $F_m$ be the algebra of subsets of $\mathbb{N}$ generated by the classes $\mathcal{B}(i, 8^m)$, $i = 0, \cdots, 8^m - 1$. $F_m$ contains finitely many elements (subsets of $\mathcal{B}(i, 8^m)$, $i = 0, \cdots, 8^m - 1$). Since $F_m$ is a finite algebra, it is trivially a $\sigma$-algebra too.
%The algebras $F_m$, $m \geq 1$ form an increasing family, since $F_{m+1}$ contains $F_m$ as a proper subalgebra. For convenience, let us call $F =\bigcup_{m\geq 1} F_m$, $F$ is a countably infinite set (being a countable union of finite sets). Furthermore, $F$ is an algebra, but, being countably infinite cannot be a sigma-algebra.}

%\ins{For each $m\geq 1$ we have a measure $\mu^{(m)}$ on the $\sigma$-algebra $F_m$. The sequence of measures is coherent in the sense that if $A$ belongs to $F_m$, then $\mu^{(m+1)}(A) = \mu^{(m)}(A)$.
%In some sense, the previous sentence tells that the measure $\mu^{(m)}$ is the \lq\lq projection\rq\rq of $\mu^{(m)}$ on the smaller algebra $F_m$. It follows that if a set $A$ belongs to $F_{m'}$ for some $m'$, then $\mu^{(m)}(A)$ is independent of $m$ (as long as it is defined, i.e. as long as $m\geq m'$). Due to this observation, for any set $A$ in $F$ we may write $\mu(A)$ and ignore the index $m$.}
%\end{remark}

\begin{remark}
Writing the singleton set $\{n\}$ as $\bigcap_{m\geq m_0} \mathcal{B}(n,8^m)$, where $m_0$ is the smallest integer such that $n\leq 8^{m_0}$, one gets $\mu^{(\infty)}_{inv}(n)=\lim_{m\rightarrow \infty}\nu(s_0)2^{-(k+1)}2^{-(m-1)}=0$. Integers are hence associated to a trivial (zero) measure. We thus have the invariance also in the limiting case $m=\infty$, even if in a trivial form $\mu^{(\infty)}_{inv}(S^{-1}\{n\})=\mu^{(\infty)}_{inv}(\{\text{a finite set of integers}\})=0$ which equals $\mu^{(\infty)}_{inv}(\{n\})=0$. Notice however that $\mu^{(\infty)}$ is not $\sigma$--additive, and so, strictly speaking, we cannot refer to it as to a measure 
\end{remark}

In the following we will be be interested in studying the orbit of the third iterate of the Collatz map, namely $S^{\circ k}$ for $k\in\mathbb{N}$. We aim in particular at characterising the quasi-stationary dynamics, i.e. the out of equilibrium dynamics of the system prior possible absorption to the equilibrium fixed point, as identified earlier. 
The transition probabilities $q_{ij}^*$ given by~\eqref{eq:Qij} have been computed using $S$, i.e. the first iterate of the map. A natural question that arise is how these latter quantities relate to the transition probabilities of the stochastic process defined through $S^{\circ k}$, for a generic $k\in\mathbb{N}$. This is established in the following proposition.  

\begin{proposition}
 \label{prop:invariance}
Let
\begin{equation*}
q^{(k)}_{ij}:=\frac{\mu_{inv}\left[\mathcal{B}(i,8)\cap S^{\circ -k}\mathcal{B}(j,8)\right]}{\mu_{inv}\left[\mathcal{B}(i,8)\right]}\, (i,j=0,\dots,7)\, ,
\end{equation*}
be the transition probability of the Markov process built from $S^{\circ k}$, for $k\in\mathbb{N}$. Then
\begin{equation}
\label{eq:invariance}
q^{(k)}_{ij}=[(q^{*})^k]_{ij}\, (i,j=0,\dots,7)\, .
\end{equation}
\end{proposition}

\proof
The proof is straightforward and relies on the Chapman-Kolmogorov equation. Let $p_{k,h}(j|i)$ be the probability to be in $\mathcal{B}(j,8)$ after $k$ steps, assuming that the system is  in $\mathcal{B}(i,8)$ after $0\leq h<k$ steps. Then
\begin{equation*}
p_{k, h}(j|i)=\sum_{m_{h+1},\dots,m_{k-1}}p_{k, k-1}(j|m_{k-1})p_{k-1, k-2}(m_{k-1}|m_{k-2})\dots p_{h+1, h}(m_{h+1}|i)\, .
\end{equation*}
Let us observe that
\begin{equation*}
p_{l+1, l}(j|i)=\frac{\mu_{inv}\left[S^{\circ -(l+1)}\mathcal{B}(j,8)\cap S^{\circ -l}\mathcal{B}(i,8)\right]}{\mu_{inv}\left[S^{\circ -l}\mathcal{B}(i,8)\right]}=\frac{\mu_{inv}\left[S^{\circ -l}\left(S^{\circ -1}\mathcal{B}(j,8)\cap \mathcal{B}(i,8)\right)\right]}{\mu_{inv}\left[S^{\circ -l}\mathcal{B}(i,8)\right]}\, ,
\end{equation*}
hence using the invariance
\begin{equation*}
p_{l+1, l}(j|i)=\frac{\mu_{inv}\left[S^{\circ -1}\mathcal{B}(j,8)\cap \mathcal{B}(i,8)\right]}{\mu_{inv}\left[\mathcal{B}(i,8)\right]}=q^*_{ij}\, .
\end{equation*}

Because $p_{k, 0}(j|i)=q^{(k)}_{ij}$ we obtain:
\begin{equation*}
q^{(k)}_{ij}=\sum_{m_1,\dots,m_{k-1}}q^*_{m_{k-1} j} q^{*}_{m_{k-2} m_{k-1}}\dots q^*_{i m_{1}}=[(q^*)^k]_{ij}\, .
\end{equation*}
\endproof

The consequences of the above Proposition are of paramount importance for the forthcoming discussion. In particular, it is possible to draw conclusion of the dynamics of $S$ (and hence on its restriction in $mod 8$) by iterating the Markov process defined by the transition probabilities $q^*_{ij}$. We recall once again that we are here interested in shedding light onto the out-of-equilibrium dynamics of  $S$ before the deterministic trajectories reach their asymptotic attractor. As we shall see, by exploiting Proposition~\ref{prop:invariance} it will be possible to deduce effective constraints that are to be matched by the deterministic application.

\subsection{Computing the transition probabilities.}
\label{ssec:transprob}

We are now in a position to explicitly compute the transition probabilities $q^*_{ij}$, formally given by Eq.~\eqref{eq:Qij}. The latter transition probabilities define the entries of the transition matrix $Q^*$ which is entirely specified by the following Proposition. 

\begin{proposition}
 \label{prop:markov}
 The matrix $Q^*$ is stochastic. The entries $q^*_{ij}$, with $i$ and $j$ in $\{0,\dots,7\}$, determine the probability to reach the target class $j$, from class $i$ and read: 
  \begin{eqnarray}
q^*_{0j}&=&\frac{1}{8}, \forall j\in\{0,\dots,7\}  \label{eq:Pk1k20}\\
q^*_{1j}&=&\frac{1}{4}, \forall j\in\{1,3,5,7\}\text{ and $0$ otherwise}\label{eq:Pk1k21}\\
q^*_{2j}&=&\frac{1}{4}, \forall j\in\{0,2,4,6\}\text{ and $0$ otherwise}\label{eq:Pk1k22}\\
q^*_{3j}&=&\frac{1}{2}, \forall j\in\{0,4\}\text{ and $0$ otherwise}\label{eq:Pk1k23}\\
q^*_{4j}&=&\frac{1}{4}, \forall j\in\{0,2,4,6\}\text{ and $0$ otherwise}\label{eq:Pk1k24}\\
q^*_{5j}&=&\frac{1}{4}, \forall j\in\{0,2,4,6\}\label{eq:Pk1k25}\\
q^*_{6j}&=&\frac{1}{4}, \forall j\in\{1,3,5,7\}\text{ and $0$ otherwise}\label{eq:Pk1k26}\\
q^*_{7j}&=&\frac{1}{2}, \forall j\in\{2,6\}\text{ and $0$ otherwise} \label{eq:Pk1k27}\, ,
\end{eqnarray}
or, equivalently,  in matrix notation:
\begin{equation}
 \label{eq:matrixQstar}
 Q^*=\left(
\begin{matrix}
 \frac{1}{8} & \frac{1}{8} & \frac{1}{8} & \frac{1}{8} & \frac{1}{8} & \frac{1}{8} & \frac{1}{8} & \frac{1}{8}\\
0 & \frac{1}{4} & 0 & \frac{1}{4} &0 & \frac{1}{4} &0 & \frac{1}{4}\\
\frac{1}{4} & 0 & \frac{1}{4} &0 & \frac{1}{4} &0 & \frac{1}{4} & 0\\
\frac{1}{2} & 0 & 0 & 0 &\frac{1}{2} & 0 & 0 & 0\\
\frac{1}{4} & 0 & \frac{1}{4} &0 & \frac{1}{4} &0 & \frac{1}{4} & 0\\
\frac{1}{4} & 0 & \frac{1}{4} &0 & \frac{1}{4} &0 & \frac{1}{4} & 0\\
0 & \frac{1}{4} & 0 & \frac{1}{4} &0 & \frac{1}{4} &0 & \frac{1}{4}\\
0 & 0 & \frac{1}{2} & 0 &0& 0 & \frac{1}{2} & 0\\
\end{matrix}
\right)\, .
\end{equation}

 \end{proposition}

\proof
Let us begin by proving that $Q^*$ is a stochastic matrix, that is $\sum_jq^*_{ij}=1$ for all $i=0,\dots,7$. Recall that classes $\mathcal{B}(l,64)$ are disjoint and their union defines the whole set of integer numbers. We hence get from Eqs.~\eqref{eq:congclass864}
\begin{equation*}
\bigcup_{j=0}^7\left(\mathcal{B}(i,8)\cap S^{-1}\mathcal{B}(j,8)\right)=\mathcal{B}(i,8)\cap \left(\bigcup_{j=0}^7 S^{-1}\mathcal{B}(j,8)\right)=\mathcal{B}(i,8)\cap \left(\bigcup_{l=0}^{63} \mathcal{B}(l,64)\right)=\mathcal{B}(i,8)\cap \mathbb{N}=\mathcal{B}(i,8)\, ,
\end{equation*}
and thus from Eq.~\eqref{eq:Qij}:
\begin{equation*}
\sum_{j=0}^7 q^*_{ij}=\sum_{j=0}^7\frac{\mu_{inv}\left[\mathcal{B}(i,8)\cap S^{-1}\mathcal{B}(j,8)\right]}{\mu_{inv}\left[\mathcal{B}(i,8)\right]}=\frac{1}{\mu_{inv}\left[\mathcal{B}(i,8)\right]}\sum_{j=0}^7\mu_{inv}\left[\mathcal{B}(i,8)\cap S^{-1}\mathcal{B}(j,8)\right]=\frac{1}{\mu_{inv}\left[\mathcal{B}(i,8)\right]}\mu_{inv}\left[\mathcal{B}(i,8)\right]=1\, .
\end{equation*}

The remaining part of the proof is straightforward.  To compute the values of $q^*_{ij}$ given by Eq.~\eqref{eq:Qij}, one makes use of Proposition~\ref{prop:smenouno}. Moreover, we recall that $\mu_{inv}[\mathcal{B}(i,8)]=1/6$ if $i=0,2,4,6$ and $1/12$ if $i=1,3,5,7$ and $\mu_{inv}[\mathcal{B}(i,64)]=\mu_{inv}[\mathcal{B}(j,8)]/8$ if $i\equiv j \mod 8$,  as follows the invariant probability measure $\mu_{inv}$. 

Consider $q^*_{0j}$:
\begin{eqnarray*}
q^*_{0j}&=&\frac{\mu_{inv}\left[\mathcal{B}(0,8)\cap S^{-1}\mathcal{B}(j,8)\right]}{\mu_{inv}\left[\mathcal{B}(0,8)\right]}=\frac{1/6\times 1/8}{1/6}=\frac{1}{8}\quad (j=0,\dots,7)\, .
\end{eqnarray*}

Let us now turn to considering $q^*_{1j}$:
\begin{eqnarray*}
q^*_{10}&=&\frac{\mu_{inv}\left[\mathcal{B}(1,8)\cap S^{-1}\mathcal{B}(0,8)\right]}{\mu_{inv}\left[\mathcal{B}(1,8)\right]}=0\notag\\
q^*_{11}&=&\frac{\mu_{inv}\left[\mathcal{B}(1,8)\cap S^{-1}\mathcal{B}(1,8)\right]}{\mu_{inv}\left[\mathcal{B}(1,8)\right]}=\frac{\mu_{inv}\left[\mathcal{B}(1,64)\cup \mathcal{B}(33,64)\right]}{\mu_{inv}\left[\mathcal{B}(1,8)\right]}=\frac{2\times 1/12\times 1/8}{1/12}=\frac{1}{4}\notag\\
q^*_{12}&=&\frac{\mu_{inv}\left[\mathcal{B}(1,8)\cap S^{-1}\mathcal{B}(2,8)\right]}{\mu_{inv}\left[\mathcal{B}(1,8)\right]}=0\notag\\
q^*_{13}&=&\frac{\mu_{inv}\left[\mathcal{B}(1,8)\cap S^{-1}\mathcal{B}(3,8)\right]}{\mu_{inv}\left[\mathcal{B}(1,8)\right]}=\frac{\mu_{inv}\left[\mathcal{B}(25,64)\cup \mathcal{B}(57,64)\right]}{\mu_{inv}\left[\mathcal{B}(1,8)\right]}=\frac{2\times 1/12\times 1/8}{1/12}=\frac{1}{4}\notag\\
q^*_{14}&=&\frac{\mu_{inv}\left[\mathcal{B}(1,8)\cap S^{-1}\mathcal{B}(4,8)\right]}{\mu_{inv}\left[\mathcal{B}(1,8)\right]}=0\notag\\
q^*_{15}&=&\frac{\mu_{inv}\left[\mathcal{B}(1,8)\cap S^{-1}\mathcal{B}(5,8)\right]}{\mu_{inv}\left[\mathcal{B}(1,8)\right]}=\frac{\mu_{inv}\left[\mathcal{B}(17,64)\cup \mathcal{B}(49,64)\right]}{\mu_{inv}\left[\mathcal{B}(1,8)\right]}=\frac{2\times 1/12\times 1/8}{1/12}=\frac{1}{4}\notag\\
q^*_{16}&=&\frac{\mu_{inv}\left[\mathcal{B}(1,8)\cap S^{-1}\mathcal{B}(6,8)\right]}{\mu_{inv}\left[\mathcal{B}(1,8)\right]}=0\notag\\
q^*_{17}&=&\frac{\mu_{inv}\left[\mathcal{B}(1,8)\cap S^{-1}\mathcal{B}(7,8)\right]}{\mu_{inv}\left[\mathcal{B}(1,8)\right]}=\frac{\mu_{inv}\left[\mathcal{B}(9,64)\cup \mathcal{B}(41,64)\right]}{\mu_{inv}\left[\mathcal{B}(1,8)\right]}=\frac{2\times 1/12\times 1/8}{1/12}=\frac{1}{4}\, .
\end{eqnarray*}

Similarly for $q^*_{2j}$ we have:
\begin{eqnarray*}
q^*_{20}&=&\frac{\mu_{inv}\left[ \mathcal{B}(2,8)\cap S^{-1} \mathcal{B}(0,8)\right]}{\mu_{inv}\left[ \mathcal{B}(2,8)\right]}=\frac{\mu_{inv}\left[\mathcal{B}(10,64)\cup \mathcal{B}(42,64)\right]}{\mu_{inv}\left[ \mathcal{B}(2,8)\right]}=\frac{2\times 1/6\times 1/8}{1/6}=\frac{1}{4}\notag\\
q^*_{21}&=&\frac{\mu_{inv}\left[ \mathcal{B}(2,8)\cap S^{-1} \mathcal{B}(1,8)\right]}{\mu_{inv}\left[ \mathcal{B}(2,8)\right]}=0\notag\\
q^*_{22}&=&\frac{\mu_{inv}\left[ \mathcal{B}(2,8)\cap S^{-1} \mathcal{B}(2,8)\right]}{\mu_{inv}\left[ \mathcal{B}(2,8)\right]}=\frac{\mu_{inv}\left[\mathcal{B} (2,64)\cup \mathcal{B}(34,64)\right]}{\mu_{inv}\left[ \mathcal{B}(2,8)\right]}=\frac{2\times 1/6\times 1/8}{1/6}=\frac{1}{4}\notag\\
q^*_{23}&=&\frac{\mu_{inv}\left[ \mathcal{B}(2,8)\cap S^{-1} \mathcal{B}(3,8)\right]}{\mu_{inv}\left[ \mathcal{B}(2,8)\right]}=0\notag\\
q^*_{24}&=&\frac{\mu_{inv}\left[ \mathcal{B}(2,8)\cap S^{-1} \mathcal{B}(4,8)\right]}{\mu_{inv}\left[ \mathcal{B}(2,8)\right]}=\frac{\mu_{inv}\left[\mathcal{B}(26,64)\cup \mathcal{B}(58,64)\right]}{\mu_{inv}\left[ \mathcal{B}(2,8)\right]}=\frac{2\times 1/6\times 1/8}{1/6}=\frac{1}{4}\notag\\
q^*_{25}&=&\frac{\mu_{inv}\left[ \mathcal{B}(2,8)\cap S^{-1} \mathcal{B}(5,8)\right]}{\mu_{inv}\left[ \mathcal{B}(2,8)\right]}=0\notag\\
q^*_{26}&=&\frac{\mu_{inv}\left[ \mathcal{B}(2,8)\cap S^{-1} \mathcal{B}(6,8)\right]}{\mu_{inv}\left[ \mathcal{B}(2,8)\right]}=\frac{\mu_{inv}\left[\mathcal{B}(18,64)\cup \mathcal{B}(50,64)\right]}{\mu_{inv}\left[ \mathcal{B}(2,8)\right]}=\frac{2\times 1/6\times 1/8}{1/6}=\frac{1}{4}\notag\\
q^*_{27}&=&\frac{\mu_{inv}\left[ \mathcal{B}(2,8)\cap S^{-1} \mathcal{B}(7,8)\right]}{\mu_{inv}\left[ \mathcal{B}(2,8)\right]}=0\, .
\end{eqnarray*}

Then for  $q^*_{3j}$ we get:
\begin{eqnarray*}
q^*_{30}&=&\frac{\mu_{inv}\left[ \mathcal{B}(3,8)\cap S^{-1} \mathcal{B}(0,8)\right]}{\mu_{inv}\left[ \mathcal{B}(3,8)\right]}=\frac{\mu_{inv}\left[\mathcal{B}(3,64)\cup \mathcal{B}(19,64)\cup\mathcal{B}(35,64)\cup \mathcal{B}(56,64)\right]}{\mu_{inv}\left[ \mathcal{B}(3,8)\right]}=\frac{4\times 1/12\times 1/8}{1/12}=\frac{1}{2}\notag\\
q^*_{31}&=&\frac{\mu_{inv}\left[ \mathcal{B}(3,8)\cap S^{-1} \mathcal{B}(1,8)\right]}{\mu_{inv}\left[ \mathcal{B}(3,8)\right]}=0\notag\\
q^*_{32}&=&\frac{\mu_{inv}\left[ \mathcal{B}(3,8)\cap S^{-1} \mathcal{B}(2,8)\right]}{\mu_{inv}\left[ \mathcal{B}(3,8)\right]}=0\notag\\
q^*_{33}&=&\frac{\mu_{inv}\left[ \mathcal{B}(3,8)\cap S^{-1} \mathcal{B}(3,8)\right]}{\mu_{inv}\left[ \mathcal{B}(3,8)\right]}=0\notag\\
q^*_{34}&=&\frac{\mu_{inv}\left[ \mathcal{B}(3,8)\cap S^{-1} \mathcal{B}(4,8)\right]}{\mu_{inv}\left[ \mathcal{B}(3,8)\right]}=\frac{\mu_{inv}\left[\mathcal{B}(11,64)\cup \mathcal{B}(27,64)\cup\mathcal{B}(43,64)\cup \mathcal{B}(59,64)\right]}{\mu_{inv}\left[ \mathcal{B}(3,8)\right]}=\frac{4\times 1/12\times 1/8}{1/12}=\frac{1}{2}\notag\\
q^*_{35}&=&\frac{\mu_{inv}\left[ \mathcal{B}(3,8)\cap S^{-1} \mathcal{B}(5,8)\right]}{\mu_{inv}\left[ \mathcal{B}(3,8)\right]}=0\notag\\
q^*_{36}&=&\frac{\mu_{inv}\left[ \mathcal{B}(3,8)\cap S^{-1} \mathcal{B}(6,8)\right]}{\mu_{inv}\left[ \mathcal{B}(3,8)\right]}=0\notag\\
q^*_{37}&=&\frac{\mu_{inv}\left[ \mathcal{B}(3,8)\cap S^{-1} \mathcal{B}(7,8)\right]}{\mu_{inv}\left[ \mathcal{B}(3,8)\right]}=0\, ,
\end{eqnarray*}

For $q^*_{4j}$:
\begin{eqnarray*}
q^*_{40}&=&\frac{\mu_{inv}\left[ \mathcal{B}(4,8)\cap S^{-1} \mathcal{B}(0,8)\right]}{\mu_{inv}\left[ \mathcal{B}(4,8)\right]}=\frac{\mu_{inv}\left[\mathcal{B}(52,64)\cup \mathcal{B}(20,64)\right]}{\mu_{inv}\left[ \mathcal{B}(4,8)\right]}=\frac{2\times 1/6\times 1/8}{1/6}=\frac{1}{4}\notag\\
q^*_{41}&=&\frac{\mu_{inv}\left[ \mathcal{B}(4,8)\cap S^{-1} \mathcal{B}(1,8)\right]}{\mu_{inv}\left[ \mathcal{B}(4,8)\right]}=0\notag\\
q^*_{42}&=&\frac{\mu_{inv}\left[ \mathcal{B}(4,8)\cap S^{-1} \mathcal{B}(2,8)\right]}{\mu_{inv}\left[ \mathcal{B}(4,8)\right]}=\frac{\mu_{inv}\left[\mathcal{B}(12,64)\cup \mathcal{B}(44,64)\right]}{\mu_{inv}\left[ \mathcal{B}(4,8)\right]}=\frac{2\times 1/6\times 1/8}{1/6}=\frac{1}{4}\notag\\
q^*_{43}&=&\frac{\mu_{inv}\left[ \mathcal{B}(4,8)\cap S^{-1} \mathcal{B}(3,8)\right]}{\mu_{inv}\left[ \mathcal{B}(4,8)\right]}=0\notag\\
q^*_{44}&=&\frac{\mu_{inv}\left[ \mathcal{B}(4,8)\cap S^{-1} \mathcal{B}(4,8)\right]}{\mu_{inv}\left[ \mathcal{B}(4,8)\right]}=\frac{\mu_{inv}\left[\mathcal{B} (4,64)\cup \mathcal{B}(36,64)\right]}{\mu_{inv}\left[ \mathcal{B}(4,8)\right]}=\frac{2\times 1/6\times 1/8}{1/6}=\frac{1}{4}\notag\\
q^*_{45}&=&\frac{\mu_{inv}\left[ \mathcal{B}(4,8)\cap S^{-1} \mathcal{B}(5,8)\right]}{\mu_{inv}\left[ \mathcal{B}(4,8)\right]}=0\notag\\
q^*_{46}&=&\frac{\mu_{inv}\left[ \mathcal{B}(4,8)\cap S^{-1} \mathcal{B}(6,8)\right]}{\mu_{inv}\left[ \mathcal{B}(4,8)\right]}=\frac{\mu_{inv}\left[\mathcal{B}(28,64)\cup \mathcal{B}(60,64)\right]}{\mu_{inv}\left[ \mathcal{B}(4,8)\right]}=\frac{2\times 1/6\times 1/8}{1/6}=\frac{1}{4}\notag\\
q^*_{47}&=&\frac{\mu_{inv}\left[ \mathcal{B}(4,8)\cap S^{-1} \mathcal{B}(7,8)\right]}{\mu_{inv}\left[ \mathcal{B}(4,8)\right]}=0\, .
\end{eqnarray*}

For $q^*_{5j}$:
\begin{eqnarray*}
q^*_{50}&=&\frac{\mu_{inv}\left[ \mathcal{B}(5,8)\cap S^{-1} \mathcal{B}(0,8)\right]}{\mu_{inv}\left[ \mathcal{B}(5,8)\right]}=\frac{\mu_{inv}\left[\mathcal{B}(21,64)\cup \mathcal{B}(53,64)\right]}{\mu_{inv}\left[ \mathcal{B}(5,8)\right]}=\frac{2\times 1/12\times 1/8}{1/12}=\frac{1}{4}\notag\\
q^*_{51}&=&\frac{\mu_{inv}\left[ \mathcal{B}(5,8)\cap S^{-1} \mathcal{B}(1,8)\right]}{\mu_{inv}\left[ \mathcal{B}(5,8)\right]}=0\notag\\
q^*_{52}&=&\frac{\mu_{inv}\left[ \mathcal{B}(5,8)\cap S^{-1} \mathcal{B}(2,8)\right]}{\mu_{inv}\left[ \mathcal{B}(5,8)\right]}=\frac{\mu_{inv}\left[\mathcal{B}(13,64)\cup \mathcal{B}(45,64)\right]}{\mu_{inv}\left[ \mathcal{B}(5,8)\right]}=\frac{2\times 1/12\times 1/8}{1/12}=\frac{1}{4}\notag\\
q^*_{53}&=&\frac{\mu_{inv}\left[ \mathcal{B}(5,8)\cap S^{-1} \mathcal{B}(3,8)\right]}{\mu_{inv}\left[ \mathcal{B}(5,8)\right]}=0\notag\\
q^*_{54}&=&\frac{\mu_{inv}\left[ \mathcal{B}(5,8)\cap S^{-1} \mathcal{B}(4,8)\right]}{\mu_{inv}\left[ \mathcal{B}(5,8)\right]}=\frac{\mu_{inv}\left[\mathcal{B} (5,64)\cup \mathcal{B}(37,64)\right]}{\mu_{inv}\left[ \mathcal{B}(5,8)\right]}=\frac{2\times 1/12\times 1/8}{1/12}=\frac{1}{4}\notag\\
q^*_{55}&=&\frac{\mu_{inv}\left[ \mathcal{B}(5,8)\cap S^{-1} \mathcal{B}(5,8)\right]}{\mu_{inv}\left[ \mathcal{B}(5,8)\right]}=0\notag\\
q^*_{56}&=&\frac{\mu_{inv}\left[ \mathcal{B}(5,8)\cap S^{-1} \mathcal{B}(6,8)\right]}{\mu_{inv}\left[ \mathcal{B}(5,8)\right]}=\frac{\mu_{inv}\left[\mathcal{B}(29,64)\cup \mathcal{B}(61,64)\right]}{\mu_{inv}\left[ \mathcal{B}(5,8)\right]}=\frac{2\times 1/12\times 1/8}{1/12}=\frac{1}{4}\notag\\
q^*_{57}&=&\frac{\mu_{inv}\left[ \mathcal{B}(5,8)\cap S^{-1} \mathcal{B}(7,8)\right]}{\mu_{inv}\left[ \mathcal{B}(5,8)\right]}=0\, .
\end{eqnarray*}

For $q^*_{6j}$:
\begin{eqnarray*}
q^*_{60}&=&\frac{\mu_{inv}\left[\mathcal{B}(6,8)\cap S^{-1}\mathcal{B}(0,8)\right]}{\mu_{inv}\left[\mathcal{B}(6,8)\right]}=0\notag\\
q^*_{61}&=&\frac{\mu_{inv}\left[\mathcal{B}(6,8)\cap S^{-1}\mathcal{B}(1,8)\right]}{\mu_{inv}\left[\mathcal{B}(6,8)\right]}=\frac{\mu_{inv}\left[\mathcal{B}(22,64)\cup \mathcal{B}(54,64)\right]}{\mu_{inv}\left[\mathcal{B}(6,8)\right]}=\frac{2\times 1/6\times 1/8}{1/6}=\frac{1}{4}\notag\\
q^*_{62}&=&\frac{\mu_{inv}\left[\mathcal{B}(6,8)\cap S^{-1}\mathcal{B}(2,8)\right]}{\mu_{inv}\left[\mathcal{B}(6,8)\right]}=0\notag\\
q^*_{63}&=&\frac{\mu_{inv}\left[\mathcal{B}(6,8)\cap S^{-1}\mathcal{B}(3,8)\right]}{\mu_{inv}\left[\mathcal{B}(6,8)\right]}=\frac{\mu_{inv}\left[\mathcal{B}(14,64)\cup \mathcal{B}(46,64)\right]}{\mu_{inv}\left[\mathcal{B}(6,8)\right]}=\frac{2\times 1/6\times 1/8}{1/6}=\frac{1}{4}\notag\\
q^*_{64}&=&\frac{\mu_{inv}\left[\mathcal{B}(6,8)\cap S^{-1}\mathcal{B}(4,8)\right]}{\mu_{inv}\left[\mathcal{B}(6,8)\right]}=0\notag\\
q^*_{65}&=&\frac{\mu_{inv}\left[\mathcal{B}(6,8)\cap S^{-1}\mathcal{B}(5,8)\right]}{\mu_{inv}\left[\mathcal{B}(6,8)\right]}=\frac{\mu_{inv}\left[\mathcal{B}(6,64)\cup \mathcal{B}(38,64)\right]}{\mu_{inv}\left[\mathcal{B}(6,8)\right]}=\frac{2\times 1/6\times 1/8}{1/6}=\frac{1}{4}\notag\\
q^*_{66}&=&\frac{\mu_{inv}\left[\mathcal{B}(6,8)\cap S^{-1}\mathcal{B}(6,8)\right]}{\mu_{inv}\left[\mathcal{B}(6,8)\right]}=0\notag\\
q^*_{67}&=&\frac{\mu_{inv}\left[\mathcal{B}(6,8)\cap S^{-1}\mathcal{B}(7,8)\right]}{\mu_{inv}\left[\mathcal{B}(6,8)\right]}=\frac{\mu_{inv}\left[\mathcal{B}(30,64)\cup \mathcal{B}(62,64)\right]}{\mu_{inv}\left[\mathcal{B}(6,8)\right]}=\frac{2\times 1/6\times 1/8}{1/6}=\frac{1}{4}\, .
\end{eqnarray*}

And, finally, for $q^*_{7j}$:
\begin{eqnarray*}
q^*_{70}&=&\frac{\mu_{inv}\left[ \mathcal{B}(7,8)\cap S^{-1} \mathcal{B}(0,8)\right]}{\mu_{inv}\left[ \mathcal{B}(7,8)\right]}=0\notag\\
q^*_{71}&=&\frac{\mu_{inv}\left[ \mathcal{B}(7,8)\cap S^{-1} \mathcal{B}(1,8)\right]}{\mu_{inv}\left[ \mathcal{B}(7,8)\right]}=0\notag\\
q^*_{72}&=&\frac{\mu_{inv}\left[ \mathcal{B}(7,8)\cap S^{-1} \mathcal{B}(2,8)\right]}{\mu_{inv}\left[ \mathcal{B}(7,8)\right]}=\frac{\mu_{inv}\left[\mathcal{B}(7,64)\cup \mathcal{B}(23,64)\cup\mathcal{B}(39,64)\cup \mathcal{B}(55,64)\right]}{\mu_{inv}\left[ \mathcal{B}(7,8)\right]}=\frac{4\times 1/12\times 1/8}{1/12}=\frac{1}{2}\notag\\
q^*_{73}&=&\frac{\mu_{inv}\left[ \mathcal{B}(7,8)\cap S^{-1} \mathcal{B}(3,8)\right]}{\mu_{inv}\left[ \mathcal{B}(7,8)\right]}=0\notag\\
q^*_{74}&=&\frac{\mu_{inv}\left[ \mathcal{B}(7,8)\cap S^{-1} \mathcal{B}(4,8)\right]}{\mu_{inv}\left[ \mathcal{B}(7,8)\right]}=0\notag\\
q^*_{75}&=&\frac{\mu_{inv}\left[ \mathcal{B}(7,8)\cap S^{-1} \mathcal{B}(5,8)\right]}{\mu_{inv}\left[ \mathcal{B}(7,8)\right]}=0\notag\\
q^*_{76}&=&\frac{\mu_{inv}\left[ \mathcal{B}(7,8)\cap S^{-1} \mathcal{B}(6,8)\right]}{\mu_{inv}\left[ \mathcal{B}(7,8)\right]}=\frac{\mu_{inv}\left[\mathcal{B}(15,64)\cup \mathcal{B}(31,64)\cup\mathcal{B}(47,64)\cup \mathcal{B}(63,64)\right]}{\mu_{inv}\left[ \mathcal{B}(7,8)\right]}=\frac{4\times 1/12\times 1/8}{1/12}=\frac{1}{2}\notag\\
q^*_{77}&=&\frac{\mu_{inv}\left[ \mathcal{B}(7,8)\cap S^{-1} \mathcal{B}(7,8)\right]}{\mu_{inv}\left[ \mathcal{B}(7,8)\right]}=0\, ,
\end{eqnarray*}
\endproof

\begin{center}
\begin{figure}[!h]
\centering
\includegraphics[width=8cm]{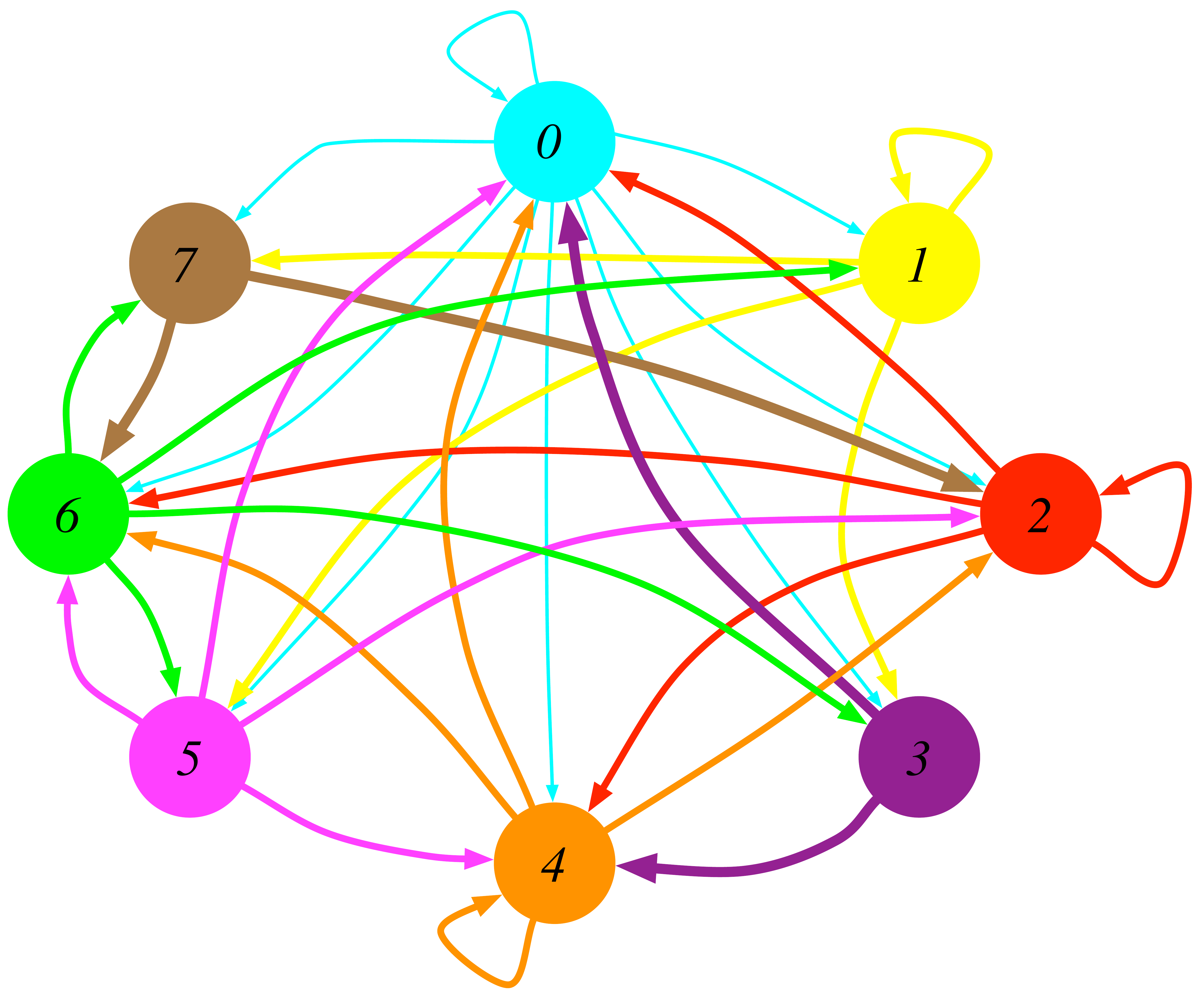} 
\caption{A schematic layout of the network that defines the allowed moves for the finite state Markov chain. The colour of the links points to the departing node while the arrows indicate the direction of the flux. The widths of the connecting arcs reflect the associated transition weights, as specified in matrix~\eqref{eq:matrixQstar}.}
\label{markovnet8}
\end{figure}
\end{center}

The finite state Markov process is graphically represented in Fig.~\ref{markovnet8}. All possible moves are shown:  the associated transition probabilities set the widths of the arcs connecting two adjacent nodes{, the thicker the arc, the larger the transition probability}.

\section{Collatz orbits are on average contracting.}
\label{sec:wCollatz}

Building on the above, we can now prove that, on average,  Collatz orbits are bounded. We begin by computing the invariant (stationary) distribution $\vec{P}_{stat}$ of the Markov process defined by the irreducible transition matrix $Q^*$. The stationary distribution $\vec{P}_{stat}$ represents a factual constraint that should be matched on average by stochastic trajectories. It will allow in turn to draw conclusions on the late time behaviour of the deterministic orbits (prior possible convergence to the Collatz cycle). We recall in fact that the transition probabilities of the Markov chain are computed by using a measure invariant (over the congruences classes of interest) under the deterministic map $S$. This allows us to study the asymptotic dynamics of the deterministic process via successive application of the Markov chain (see Proposition~\ref{prop:invariance}).

Determining the stationary distribution $\vec{P}_{stat}$ amounts to solving the following eigenvalue problem:

\begin{equation}
\label{eq:qstat}
\vec{P}_{stat}  Q^* =  \vec{P}_{stat} 
\end{equation}

In practical terms, the stationary distribution corresponds to the left eigenvector associated to the dominant eigenvalue ($\lambda=1$, due to the stochasticity constraint) of the transition matrix $Q^*$. Performing the above calculation immediately returns: 
 
\begin{equation}
 \label{eq:P_QS}
\vec{P}_{stat}=(1/6,1/12,1/6,1/12,1/6,1/12,1/6,1/12)
\end{equation}

To proceed in the analysis we prove the following Proposition:
\begin{proposition}
 \label{prop:markoverg}
 The finite states Markov process defined by the transition matrix $Q^*$ is irreducible and recurrent, hence ergodic.
\end{proposition}

\proof
The claim follows by direct inspection of $(Q^*)^2$, the second iterate of $Q^*$. In fact, $[(Q^*)^2]_{ij}\geq 1/16$ for all $i$ and $j$ in $\{0,\dots,7\}$.
\endproof

Recall now that, based on Proposition~\ref{prop:invariance}, one can elaborate on the late time dynamics of $S$ (or equivalently, its restriction in $mod8$) by virtue of the introduced Markov chain.  This observation translates in a viable strategy to prove the average contracting property of the deterministic map $S$. Roughly speaking when the map $S$ is acted on a positive integer belonging to class $\mathcal{B}(0,8)$, it drives a contraction factor equal to $1/8$. If $S$ is instead operated on natural numbers of the classes $\mathcal{B}(j,8)$, $j=1,2,4,5,6$, it results in a contraction of $3/4$. At variance,  $S(n)$ produces an expansion with rate $9/2$, for $n$ belonging to classes $\mathcal{B}(3,8)$ and $\mathcal{B}(7,8)$. Hence the obtained information on  $\vec{P}_{stat}$ allows to estimate the degree of contraction (or expansion) $f_{Q^*}$ that trajectories should, on average, produce:
\begin{equation}
 \label{eq:Qstella}
 f_{Q^*} \simeq \left(\frac{1}{8}\right)^{\omega_0}\left(\frac{3}{4}\right)^{\omega_1}\left(\frac{9}{2}\right)^{\omega_2}\, 
\end{equation}
where $\omega_0=1/6$ is the probability of being in the congruence class $j=0$, $\omega_1=3/6+2/12=2/3$ the probability of being in the congruence classes $j=1,2,4,5,6$ and $\omega_2=2/12=1/6$ the probability of being in the congruence classes $j=3,7$. Probabilities  $\omega_0,$  $\omega_1$,  $\omega_2$ follow equation  (\ref{eq:P_QS}), while the contraction/expansion factors $\frac{1}{8}$, $\frac{3}{4}$ and $\frac{9}{2}$ associated to each of congruence class are made explicit in  Eqs. (\ref{eq:the3map}).

Carrying out  the calculation yields $f_{Q^*}= 3/4<1$, thus implying in turn that the average approach to the absorbing equilibrium is contracting. Observe that this latter contracting factor is here analytically determined, at variance with previous attempts that relied on heuristic reasoning. Notice that this preliminary estimate $f_{Q^*}= 3/4$ has been obtained by just retaining the terms proportional to $n$ in the definition of $S(n)$, see Eqs. (\ref{eq:the3map}), or, equivalently, working with a sufficiently large $n$. Accounting for the constant ($n$ independent) contributions in (\ref{eq:the3map}) does not modify the conclusion that we have reached: the generic orbit is always contracting, as it is proved hereafter.

Consider in fact Eqs.~\eqref{eq:the3map} which define the map $S$ on the classes $\mathcal{B}(i,8)$. We are in particular interested in the contraction/expansion factors associated to each transition among classes. The following upper bounds can be obtained:

\begin{equation}
 \label{eq:the3map_upperbound}
 \forall n\in\mathbb{N}\quad S(n)=
 \begin{cases}
\frac{n}{8}=:c_0(n_{min}) n& \text{if $n\in \mathcal{B}(0,8)$}\\
\frac{3n+1}{4} < \frac{3}{4}(1+\frac{1}{3n_{min}})n =:c_1(n_{min}) n& \text{if $n\in \mathcal{B}(1,8)$}\\
\frac{3n+2}{4} <  \frac{3}{4}(1+\frac{2}{3n_{min}})n=:c_2(n_{min}) n & \text{if $n\in \mathcal{B}(2,8)$}\\
\frac{9n+5}{2} <  \frac{9}{2}(1+\frac{5}{9n_{min}})n=:c_3(n_{min}) n& \text{if $n\in \mathcal{B}(3,8)$}\\
\frac{3n+4}{4} <  \frac{3}{4}(1+\frac{4}{3n_{min}})n=:c_4(n_{min}) n& \text{if $n\in \mathcal{B}(4,8)$}\\
\frac{3n+1}{4}<  \frac{3}{4}(1+\frac{1}{3 n_{min}})n=:c_5(n_{min}) n& \text{if $n\in \mathcal{B}(5,8)$}\\
\frac{3n+2}{4}<  \frac{3}{4}(1+\frac{2}{3n_{min}})n=:c_6(n_{min}) n& \text{if $n\in \mathcal{B}(6,8)$}\\
\frac{9n+5}{2}<  \frac{9}{2}(1+\frac{5}{9n_{min}})n=:c_7(n_{min}) n& \text{if $n\in \mathcal{B}(7,8)$\, ,}
 \end{cases}
\end{equation}
where the constants $c_i(n_{min})$ are defined by the rightmost hand sides of the previous equations. Here $n_{min}$ stands for the smallest integer visited by the system in its quasi-stationary state (i.e. before it eventually hits the absorbing Collatz cycle, if this is the case). Since  $\{1,2,4\}$ belong to the Collatz cycle, and because we are solely focusing on the dynamics that precedes the possible convergence to the Collatz cycle, we will set  $n_{min}=3$~\footnote{In principle one could also consider the additional constraint on $n_{min}$, imposed by the class of  relative pertinence.}. With this choice, one gets:
\begin{equation}
 \label{eq:thecis}
c_0=\frac{1}{8}\, , c_1=\frac{3}{4}\frac{10}{9}\, ,c_2=\frac{3}{4}\frac{11}{9}\, , c_3=\frac{9}{2}\frac{32}{27}\, ,c_4=\frac{3}{4}\frac{13}{9}\, ,c_5=\frac{3}{4}\frac{10}{9}\, ,c_6=\frac{3}{4}\frac{11}{9}\text{ and }c_7=\frac{9}{2}\frac{32}{27}\, ,
\end{equation}
where the explicit reference to  $n_{min}=3$ has been dropped in the definition of the symbols $c_i$, $i=0,...,7$. 

We have therefore:
\begin{equation}
 \label{eq:Qstella1}
 f_{Q^*} \le c_0^{\Omega_0} c_1^{\Omega_1} c_2^{\Omega_2} c_3^{\Omega_3}c_4^{\Omega_4}\, ,
 \end{equation}
where  $\Omega_0=\omega_0=1/6$ is the probability of being in the congruence class $j=0$, $\Omega_1=1/6$ the probability of being in the congruence classes $j=1,5$, 
$\Omega_2=1/3$ the probability of being in the congruence classes $j=2,6$, $\Omega_3=1/6$ the probability of being in the congruence classes $j=3,7$ and $\Omega_4=1/6$ the probability of being in the congruence class $j=4$. Performing the calculation yields $f_{Q^*} \le 0.8926$. The dynamics of $S$ is therefore contracting and trajectories are on average attracted towards the three fixed points as identified above, namely the entries of the Collatz cycle $\{1,2,4\}$.

The previous observation can be also made rigorous  at the level of single trajectories.
Consider in fact a generic Markov chain with state space $X$ and transition matrix $Q$. Assume the Markov chain to be irreducible and positive recurrent. 
Let $\pi$ denote the unique invariant probability measure and consider
a non-negative function $f : X \mapsto \mathbb{R}$, summable with respect to $\pi$. Then for a.e $n_0 \in X$, we recall that 
\begin{equation}
\lim_{k\rightarrow\infty} \frac{1}{k}\sum_{j=0}^{k-1}f(S^{\circ j}(n_0)) =\sum f(x) \pi(x)\, .
\end{equation}

Taking $f=\delta_i$, $i\in\{0,\dots,7\}$, the latter result implies that the entries of the stationary distribution represent the fraction of time spent by the Markov chain in each of the eight classes, more precisely for almost every orbit of the stochastic process, $X_j$, we have for all $i\in\{0,\dots,7\}$:
\begin{equation*}
\lim_{k\rightarrow \infty}\frac{\#\{0\leq j \leq k-1: X_j=i\}}{k}=\lim_{k\rightarrow \infty}\frac{1}{k}\sum_{j=0}^{k-1}\delta_{X_j=i}= (P_{stat})_i\, .
\end{equation*}

Equipped with the above one can prove that almost all orbits generated by $S$ are bound to contract and hence converge to the Collatz cycle.
To this end let us assume the existence of an initial datum $n_0 \in  \mathbb{N}$ associated to a diverging Collatz orbit, $n_k=S^{\circ k}(n_0)$, $k\geq 0$. 
From definition  (\ref{eq:thecis}) it obviously follows $S(n)<c_i n$ for $n\in\mathcal{B}(i,8)$ and thus
\begin{equation}
 \label{eq:nkcis}
n_k=S(n_{k-1})<c_{i_{k-1}}n_{k-1}=c_{i_{k-1}}S(n_{k-2})<c_{i_{k-1}}c_{i_{k-2}}n_{k-2}<\dots< c_{i_{k-1}}\dots c_{i_{0}}n_0\, .
\end{equation}

We now proceed by defining the quantity $Y_k=\log c_{i_k}$. As the trajectory stemming from $n_0$ is assumed by hypothesis diverging, and hence constituted by an infinite number of entries, one gets:

\begin{equation}
\label{eq:thmerg}
 \lim_{k\rightarrow \infty}\frac{1}{k}\sum_{j=0}^{k-1}Y_l=\frac{\log
 c_0}{6}+\frac{\log c_1}{12}+\frac{\log c_2}{6}+\frac{\log
 c_3}{12}+\frac{\log c_4}{6}+\frac{\log c_5}{12}+\frac{\log
 c_6}{6}+\frac{\log c_7}{12}\sim -0.1136=\alpha\, .
\end{equation}
where the factors $1/6$ and $1/12$ that weights $c_i$ follows from the stationary distribution $\vec{P}_{stat}$ as computed above.

It is therefore always possible to choose an integer $k(n_0)>0$ such that, for all $k>k(n_0)$ one, has
\begin{equation*}
\sum_{j=0}^{k-1}Y_l \leq \frac{\alpha}{2}k\, ,
\end{equation*}

Notice that in the above bound the factor $1/2$ is arbitrary and does not bear any degree of specificity.  It hence follows:

\begin{equation*}
c_{i_{k-1}}\dots c_{i_{0}} =e^{\sum_{j=0}^{k-1}Y_l}\leq \beta^k \quad \forall k>k(n_0)\, ,
\end{equation*}
where $\beta=e^{\alpha/2}\sim 0.944$.

Summing up, under the assumption of the existence of a {\it diverging orbit} we have proved that, for a sufficiently large {\it stopping time}~\cite{Lagarias} $k$, one has: 
\begin{equation*}
S^{\circ k}(n_0)=n_k \leq \beta^k n_0\, .
\end{equation*}
Since $\beta<1$, this contradicts the assumption of dealing with a diverging orbit. The above result holds for $\mu_{inv}$--almost every initial conditions $n_0$. In other words,   for $\mu_{inv}$--almost every initial conditions $n_0$ no diverging orbits can exist. 

The above results follow a dynamical constraint on the equilibrium of $S$ as obtained by partitioning the integers in $8$ congruence classes. The visiting frequency are in fact obtained by computing the stationary distribution of a Markov analogue of the deterministic dynamics which runs on a finite alphabet of $8$ states.  What is going to happen if the analysis is progressively refined to the smaller scales, by working on the classes $\mathcal{B}(i,8^m)$, for any given choice of $m$? We shall adapt the Markov approach to account for this generalisation and prove that orbits are contracting, at any given degree resolution, i.e. when sampling the equilibrium on the equivalence classes $\mathcal{B}(i,8^m)$ for $m$ large as sought. Remarkably, the computed upper bound for the contraction factor is independent on $m$ and equal to the value obtained in eq.~\eqref{eq:Qstella1}.

\section{The Collatz dynamics is contracting at the finest scale.}
\label{sec:wCollatzfinescale}

In this section we will prove that Collatz orbits are on average contracting when seen on the equivalence classes $\mathcal{B}(i,8^m)$, $i\in\{0,\dots, 8^m-1\}$ for any $m>1$ . Moreover, the computed upper bound for the contracting factors is identical to that obtained above, when operating with classes $\mathcal{B}(i,8)$, see Eq.~\eqref{eq:Qstella1}, and thus independent of $m$.

To this aim, the first step is to compute the transition probabilities $q_{ij}(m)$. This latter quantifies the probability of reaching class $\mathcal{B}(j,8^m)$ when starting from class 
$\mathcal{B}(i,8^m)$, that is the conditional probability $P\left[S(x)\in \mathcal{B}(j,8^m)| x\in \mathcal{B}(i,8^m)\right]$. In formulae:
\begin{equation}
\label{eq:Qijm}
q_{ij}(m):=\frac{\mu_{inv}\left[\mathcal{B}(i,8^m)\cap S^{-1}\mathcal{B}(j,8^m)\right]}{\mu_{inv}\left[\mathcal{B}(i,8^m)\right]}\, (i,j=0,\dots,8^m-1)\, .
\end{equation}

Label  $Q(m)$ the $8^m \times 8^m$ matrix formed by the entries (\ref{eq:Qijm}), for any given $m$.
We can then prove the following proposition:

\begin{proposition}
 \label{prop:markovQm}
The matrices $Q(m)$ are stochastic for all $m\geq 1$. Moreover the unique stationary distribution $\vec{P}_{stat}(m)$, solution of
\begin{equation*}
\vec{P}_{stat}(m)Q(m)=\vec{P}_{stat}(m)\, ,
\end{equation*}
is the vector of $\mathbb{R}^{8^m}$ given by
\begin{equation*}
\vec{P}_{stat}(m)=(a,b,\dots,a,b)\, ,
\end{equation*}
where $a=1/6\times 1/8^{m-1}$ and $b=1/12\times 1/8^{m-1}$.
  \end{proposition}

\proof
The first claim can be proved by observing that that $S^{-1}\mathcal{B}(j,8^m)$ is made by the union of disjoint classes $\mathcal{B}(l_j,8^{m+1})$, where the set of indexes $l_j$ depends on the initial class $\mathcal{B}(j,8^m)$. Hence
\begin{eqnarray*}
\bigcup_{j=0}^{8^m-1}\left(\mathcal{B}(i,8^m)\cap S^{-1}\mathcal{B}(j,8^m)\right)&=&\mathcal{B}(i,8^m)\cap \left(\bigcup_{j=0}^{8^m-1} S^{-1}\mathcal{B}(j,8^m)\right)=\mathcal{B}(i,8^m)\cap \left(\bigcup_{l=0}^{{8^{m+1}-1}} \mathcal{B}(l,{8^{m+1}})\right)= \\
&=&\mathcal{B}(i,8^m)\cap \mathbb{N}=\mathcal{B}(i,8^m)\, ,
\end{eqnarray*}
and thus from Eq.~\eqref{eq:Qijm}:
\begin{equation*}
\sum_{j=0}^{8^m-1} q_{ij}(m)=\sum_{j=0}^{8^m-1}\frac{\mu_{inv}\left[\mathcal{B}(i,8^m)\cap S^{-1}\mathcal{B}(j,8^m)\right]}{\mu_{inv}\left[\mathcal{B}(i,8^m)\right]}=\frac{1}{\mu_{inv}\left[\mathcal{B}(i,8^m)\right]}\mu_{inv}\left[\mathcal{B}(i,8^m)\right]=1\, .
\end{equation*}

Let us now prove that $\vec{P}_{stat}(m)=(a,b,\dots,a,b)$, $a=1/6\times 1/8^{m-1}$ and $b=1/12\times 1/8^{m-1}$, is the (unique) eigenvector associated to the eigenvalue equal to unit (notice that this latter exists because $Q(m)$ is stochastic). Let us start by computing
\begin{equation}
\label{eq:PQm}
\sum_i (\vec{P}_{stat}(m))_iq_{ij}(m)=a\sum_{i\, \text{even}} q_{ij}(m)+b\sum_{i\, \text{odd}} q_{ij}(m)\, .
\end{equation}
Observe then that, apart from the normalising factor $1/8^{m-1}$, $q_{ij}(m)$ is given by the number of distinct solutions of the linear congruence equations $n\equiv i \text{ mod$8^m$}$ and $S(n)\equiv j \text{ mod$8^m$}$, namely $N_{ij}:=\#\{ \mathcal{B}(i,8^m)\cap S^{-1}\mathcal{B}(j,8^m)\}$. This allows us to rewrite Eq.~\eqref{eq:PQm} as follows
\begin{equation*}
(\vec{P}_{stat}(m)Q(m))_j =\frac{a}{8} \sum_{i\, \text{even}} N_{ij}+\frac{b}{8} \sum_{i\, \text{odd}} N_{ij}\, ,
\end{equation*}
where used has been made of the definition of the measure $\mu_{inv}^{(m)}$.

By recalling Proposition~\ref{prop:AeAo} we can now write
\begin{equation*}
\mathcal{B}(i,8^m)\cap S^{-1}\mathcal{B}(j,8^m)=\mathcal{B}(i,8^m)\cap\left(A^{(m)}_e(j)\cup A^{(m)}_o(j)\right)\, ,
\end{equation*}
where $A_e(j)$ is the union of disjoint classes $\mathcal{B}(l,8^{m+1})$ with $l$ even and $A^{(m)}_o(j)$ is the union of disjoint classes $\mathcal{B}(l,8^{m+1})$ with $l$ odd.

By invoking again Proposition~\ref{prop:AeAo} we can write
\begin{equation*}
\sum_{i\, \text{even}} N_{ij}=
\begin{cases}
5&\text{if $j$ is even}\\
3&\text{if $j$ is odd}
\end{cases}
\end{equation*}
and
\begin{equation*}
\sum_{i\, \text{odd}} N_{ij}=
\begin{cases}
6&\text{if $j$ is even}\\
2&\text{if $j$ is odd}
\end{cases}
\end{equation*}

Hence, in conclusion
\begin{equation*}
\begin{cases}
\frac{5}{8}a+\frac{6}{8}b=a & \text{if $j$ is even}\\
\frac{3}{8}a+\frac{2}{8}b=b & \text{if $j$ is odd}
\end{cases}
\end{equation*}
which returns  $a=1/6$ and $b=1/12$ as the sole non trivial solution. Observe that the factor $1/8^{m-1}$ is needed to normalise the $1$-norm of the vector to $1$.
Finally let us observe that this latter conclusion generalises  to the non relatively prime case under study the results obtained by~\cite{Matthews1,Matthews2} for  relatively prime settings.
\endproof 

\begin{remark}
As shown above, the Markov process characterised by the stochastic matrix $Q(m)$ admits a unique stationary state, with all non trivial entries. Then it is immediate to show that is is also  
irreducible and recurrent. In fact for arbitrarily large $k$,  all rows of $Q(m)^k$ are identical to the stationary distribution $\vec{P}_{stat}(m)$. Its components are therefore strictly positive. 
\end{remark}

We are now in a position to quantify the degree of contraction/expansion that characterises the deterministic dynamics, as seen on the classes $\mathcal{B}(i,8^m)$, $m\geq 1$ and $i=0,\dots,8^m-1$. This enables us  to generalise the above analysis, beyond the rather specific choice $m=1$.  To this end let us define $i^{(m)}(n)=i$ if $n\in \mathcal{B}(i,8^m)$, that is the indicator function of the classes $\mathcal{B}(i,8^m)$. Observe that the contraction/expansion $S$ is ultimately determined by $i^{(1)}(n)$ (see Eq.~\eqref{eq:the3map_upperbound}). To prove that contracting character of the map we proceed as above, and assume that the orbit stemming from $n_0$ is diverging, hence by definition made by an infinite set of entries. 

Introduce then the function $W(n)=\log c_{i^{(m)}(n)}$. Since the time spent by the system in each of the available classes is determined by the components of eigenvector with eigenvalue equal to one, for  a.e. $n_0$ we have

\begin{equation*}
 \lim_{k\rightarrow \infty}\frac{1}{k}\sum_{j=0}^{k-1}W(S^{\circ j}(n_0))=
 \left[ \frac{\sum_{k\,\in \text{\{0,2,4,6\}}}\log c_k}{6\; 8^{m-1}}+\frac{\sum_{k\, \in \text{\{1,3,5,7\}}}\log c_k}{12\; 8^{m-1}} 
 \right]8^{m-1}=\alpha\, ,
\end{equation*}

where $1/6 \times 1/8^{m-1}$ and $1/12 \times 1/8^{m-1}$ are the entries of the stationary invariant distribution $\vec{P}_{stat}(m)$ (i.e. the time spent in each of the $8^{m}$ classes)  and the overall factor $8^{m-1}$ is due to the fact that for any fixed $i^{(1)}=i$ for $i\in\{0,\dots,7\}$, there are $8^{m-1}$ possibilities to have $i^{(m)}=i \textit{ mod$8^m$}$. In conclusion the map $S$ is on average contracting on the classes $\mathcal{B}(i,8^m)$. As an additional, important remark, we notice that the upper bound for the average contracting factor as obtained when representing the dynamics of $S$ on the  $8^{m}$ equivalent classes is identical to that given by Eqs. (\ref{eq:Qstella1}). 

The remarkable conclusion is therefore that the third iterate of the Collatz map is always contracting, when seen on the equivalence classes $\mathcal{B}(i,8^m)$, for $m$ large as sought,  
and that the estimated bound for the contraction factor is independent on the classes index $m$. In other words,  we can make
 the number of classes as large as wished (and consequently reduce their size so to approach the singletons with arbitrary accuracy), while still detecting a contracting deterministic dynamics, 
 with a constant (independent on $m$.) bound for the rate of contraction. 
 As previously remarked when the limit for $m$ that goes to infinity is performed, the measure of the classes, and hence the singletons, converges to zero. Despite the fact the contracting factors stays constant for any, arbitrarily large $m$, it seems that we cannot rule out the existence of zero measure orbits that violate this constraint.
 
 \begin{remark}
In principle it would be tempting to consider a different measure that would become non trivial in the limit for $m$ that diverges to infinity. This latter measure needs however to be invariant on $S$, which implies that it should asymptotically concentrate only on the Collatz integers $1,2,4$, i.e. the supposed attractors for the deterministic dynamics. On the other hand, the invariance requirement at the coarse grained scale, necessitates dealing with a uniform (except for the weights $1/6$ and $1/12$) measure on the equivalence classes.  This latter request cannot be reconciled with the need for a non uniform measure at integers level, making impossible the search for a modified measure  that is both non trivial and invariant on the singletons. 
\end{remark}

\section{Conclusions}

In this paper we have provided an analytical argument to support the validity of the so called Collatz conjecture, a long standing problem in mathematics which dates back to 1937. The analysis builds on three main pillars. In short, we (i) introduced the (forward) third iterate of the Collatz map (so to reduce the analysis of the period 3 cycle to a search for a fixed point) and considered the equivalence classes of integer numbers modulo 8; (ii)  defined a Markov chain (based on a suitable non trivial measure) which runs on a set of finite states and whose transition probabilities reflect the deterministic map; (iii)  showed that orbits are on average contracting, as follows strict bound that combines the visiting frequencies, as derived in the framework of the aforementioned stochastic picture, and the contraction/expansion factors associated to each transition among classes. Notice that the conclusion reached holds for any level of imposed coarse graining, i.e. by computing the visiting frequencies on the partition in {\textit{mod$8^m$}} classes, with $m$ large as wished. Despite the measure introduced cannot be extended to weight individual singletons, we can proof that the Collatz dynamics is contracting on uniform partitions of the natural numbers in classes. These partitions 
can be refined to approximate singletons with suited accuracy, without eventually converging to them.

\section*{Acknowledgments}
{We would like to thank the numerous colleagues who interacted with us all along the various stages of the writing of this work. In particular we warmly thank Claudio Bonanno, Carlo Carminati, Jean-Charle Delvenne, Craig Alan Feinstein, Steffen Kionke, Shlomo Levental, Stefano Marmi, Vassilis Papanicolaou, Fran\c{c}ois Stealens and C\'edric Villani, for their insightful comments and remarks. }
The work of T.C. presents research results of the Belgian Network DYSCO (Dynamical Systems, Control, and Optimization), funded by the Interuniversity Attraction Poles Programme, initiated by the Belgian State, Science Policy Office.

\appendix
\section{Methods}

\subsection{A numerical test}

\begin{center}
\begin{figure}[!h]
\centering
\includegraphics[width=8cm]{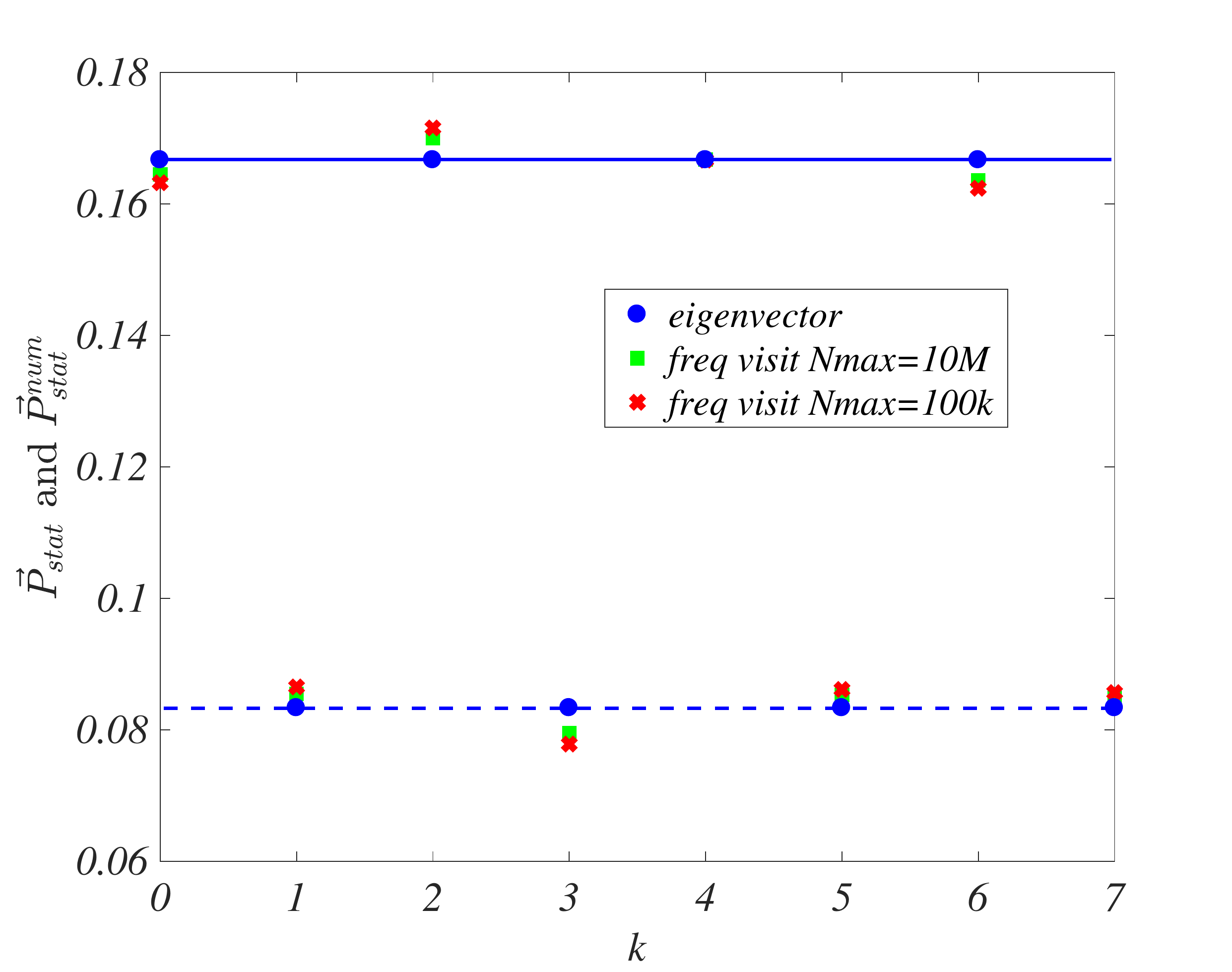}
\caption{The theoretical quasi-stationary distribution (filled blue circles) is compared to the homologous quantity obtained from a direct numerical implementation of the deterministic Collatz map. More specifically, starting from a generic integer $n$, we store the \textit{mod8} representation of the numbers obtained every $3$ consecutive iterations of $T$, before the trajectory lands on the attractive, absorbing cycle. The frequency of visits for each of the allowed classes is then reported in the Figure. Red crosses are obtained running the Collatz map for all integers in between $1$ and $n_{max}=10^5$. Green squares refer to  $n_{max}=10^7$. The blue horizontal lines are drawn as guideline for the reader, the solid line corresponds to $1/6$ and the dashed one to $1/12$. When increasing  $n_{max}$, all entries of the quasi-stationary distribution shifts consistently, although imperceptibly, towards the theoretically predicted values.  Observe that, for $n_{max} =10^7$, the largest value attained by the system is $M \simeq 6\times10^{13}$}
\label{fig2}
\end{figure}
\end{center}

The components of $\vec{P}_{stat}$ quantify the time spent by the system in each of the eight classes, in its long-lasting approach to equilibrium. To challenge this interpretative picture we compare in Figure~\ref{fig2} the prediction for $\vec{P}_{stat}$  as obtained under the stochastic picture to the quasi-stationary distribution $\vec{P}^{num}_{stat}$  recorded numerically from a direct implementation of the deterministic Collatz map. By quasi-stationary distribution we here mean the frequency of visits of the eight classes, prior absorption to the Collatz cycle. 

\end{document}